\documentclass[a4paper, reqno,11pt, twoside]{amsart}
\usepackage[latin1]{inputenc}
\usepackage[T1]{fontenc}
\usepackage{latexsym}
\usepackage{amsmath,amsthm,amssymb,amsfonts}
\usepackage{enumerate}
\usepackage{graphicx}

\usepackage[arrow,matrix,tips,curve]{xy}

\theoremstyle{definition}
\newtheorem{definition}{Definition}[section]

\newtheorem{example}[definition]{Example}

\theoremstyle{remark}

\newtheorem{remark}[definition]{Remark}

\theoremstyle{plain}
\newtheorem{lemma}[definition]{Lemma}

\newtheorem{theorem}[definition]{Theorem}

\newtheorem{corr}[definition]{Corollary}

\newcommand{\M}{\mathcal{M}}
\newcommand{\rest}[1]{\rho_{{#1}}}

\newcommand{\UNI}{\mathcal{A}}
\newcommand{\divides}[2]{{#1 \left\lvert {#2} \right.}}

\newcommand{\ddivides}[2]{{#1 \left\lvert\lvert {#2} \right.}}
\newcommand{\sqbij}{\Phi}

\newcommand{\SQA}{\mathfrak{D}}
\newcommand{\GTRUNC}{\mathfrak{S}}

\newcommand{\norm}[1]{{ \left\lvert {#1} \right\rvert}}
\newcommand{\tdeg}[1]{\norm{#1}}
\newcommand{\ord}{\mathrm{ord}}

\newcommand{\supp}{\mathrm{supp}}
\newcommand{\oksets}{\mathcal{O}}
\newcommand{\upsets}{\mathcal{U}}
\newcommand{\Tor}{\mathrm{Tor}}

\newcommand{\pr}{\mathrm{pr}}
\newcommand{\sqfree}{\tilde{\Nat}}
\newcommand{\finsubs}{\mathcal{FIN}}

\DeclareTextSymbol{\textbackslash}{T1}{92}

\newcommand{\setC}{{\mathord{\mathbb C}}}
\newcommand{\setN}{{\mathord{\mathbb N}}}

\newcommand{\Nat}{{\mathbb{N}}}
\newcommand{\Z}{{\mathbb{Z}}}
\newcommand{\C}{{\mathbb{C}}}

\newcommand{\Primes}{{\mathbb{P}}}
\newcommand{\PP}{{\mathbb{PP}}}
\newcommand{\vektor}[1]{{\boldsymbol{#1}}}

\def\setsuchas#1#2{\left\{\,{#1}\,\vrule\,{#2}\,\right\}}
\newcommand{\set}[1]{{\{#1\}}}

\begin{document}

\title[General truncations]{The
  ring of arithmetical 
  functions with unitary convolution: General Truncations.}
\author{Jan Snellman}
\address{Department of Mathematics\\
Stockholm University\\
SE-10691 Stockholm,
Sweden}
\email{jans@matematik.su.se}




\begin{abstract}
  Let \((\UNI,+,\oplus)\) denote the ring of arithmetical functions
  with unitary convolution, and let \(V \subseteq \Nat^+\) have the
  property that for every \(v \in V\), all unitary divisors of \(v\)
  lie in \(V\). If in addition \(V\) is finite, then
  \(\UNI_V\) is an artinian
  monomial quotient of a polynomial ring in finitely many
  indeterminates, and  isomorphic to the ``Artinified''
  Stanley-Reisner ring \(\setC[\overline{\Delta(V)}]\) of a certain
  simplicial complex \(\Delta(V)\). We describe some ring-theoretical
  and homological properties of \(\UNI_V\).
\end{abstract}

\maketitle


\begin{section}{Introduction}
  The ring of arithmetical functions with Dirichlet convolution 
  \begin{displaymath}
    f * g (n) = \sum_{\divides{d}{n}} f(d) g(n/d)
  \end{displaymath}
  is
  well-studied and well understood. From a ring-theoretical point of
  view, the most important fact is the theorem by Cashwell-Everett
  \cite{Cashwell:FPS} which states that this ring is a unique
  factorization domain (it is isomorphic to the power series ring over
  \(\C\) on
  countably many variables). It is also interesting that this ring can
  be given a natural norm, with respect to which it is a normed,
  valued ring \cite{vring, dvar}. 

  The ring of arithmetical functions with unitary convolution
  \begin{displaymath}
    f * g (n) = \sum_{\substack{\divides{d}{n}\\ \gcd(d,n/d)=1}} f(d)
    g(n/d) 
  \end{displaymath}
  is, by contrast, not even a domain. However, the same norm as before
  makes it into a normed (not valued) ring \cite{Snellman:UniDivTop}
  and it is at least présimplifiable, atomic, and has bounded length
  on factorizations of a given element.

  In \cite{Snellman:Dirichlet} certain truncations of the ring of
  arithmetical functions with Dirichlet convolution was studied. These
  rings were defined as follows: fix a positive integer \(n\), and
  consider all arithmetical functions supported on
  \([n]=\set{1,2,\dots,n}\). If we modify the multiplication slightly,
  this becomes a zero-dimensional algebra, which is the monomial
  quotient of a polynomial ring (on finitely many
  indeterminates). Furthermore, the defining ideals are \emph{strongly
    stable}, so the homological properties of these truncated algebras
  are easily determined, using the Eliahou-Kervaire resolution
  \cite{Eliahou:MinRes}. A way of stating that these truncated
  algebras ``approximate'' the original algebra is to note that they
  form an inverse system (of discretely topologized Artinian
  \(\C\)-algebras) whose inverse limit is precisely the ring of
  arithmetical functions with Dirichlet convolution (and with the
  topology induced by the above mentioned norm). 
  
  Of course, the same truncation can be performed on the ring of
  arithmetical functions with unitary convolution. In this article,
  we'll be somewhat more general, and consider truncations to any
  subset \(V \subseteq \Nat^+\) that is closed w.r.t taking unitary
  divisors. Of course, the case when \(V\) is finite is the most
  interesting one. Here, the truncated algebra is again a monomial
  quotient of a polynomial ring, but the defining ideals need not be
  strongly stable. They are, however, of the form
  \(I+(x_1^2,\dots,x_v^2)\) where \(I\) is square-free, so they can be
  regarded as \emph{Artinified Stanley-Reisner rings} on certain
  simplicial complexes. As a matter of fact, we show
  (Theorem~\ref{thm:everything}) that the truncations comprise
  \emph{all} 
  Artinified Stanley-Reisner rings. 

  In a forthcoming article, we'll study the particular truncations
  \(V=[n]\), exactly as in \cite{Snellman:Dirichlet}. The defining
  ideals are now multi-stable rather than stable, so the minimal free
  resolutions are more involved. However, the associated simplicial
  complexes, as defined above, are \emph{shellable}, which makes it
  possible to give some information about the Poincaré-Betti series.
\end{section}

\begin{section}{The ring of arithmetical functions with unitary
    convolution}
  We denote the \(i\)'th prime number by \(p_i\), and the set of all
  prime numbers by \(\Primes\). The set of prime powers is denoted
  by \(\PP\). For a positive integer \(n\), we put \([n] =
  \set{1,\dots,n}\). 
  \begin{subsection}{Unitary multiplication, unitary convolution, the
      norm of an arithmetical function}
  The unitary multiplication for positive integers is defined by 
  \begin{equation}
    \label{eq:unm}
d \oplus m = 
\begin{cases}
  dm & \text{ if } \gcd(d,m)=1\\
  0 & \text{ otherwise}
\end{cases}
  \end{equation}
With this operation, \((\setN^+,\oplus)\) becomes a monoid-with-zero
(see Definition~\ref{def:monoidwzero}) . 
We write \(\ddivides{d}{n}\) (or sometimes \(d \leq_\oplus n\))
when \(d\) is a unitary divisor of
  \(n\), i.e. when \(n=d \oplus m\) for some \(m\). Then \((\setN^+,
  \leq_\oplus)\)  is a partially ordered set, and a well-order, but not
  a lattice, since for instance \(\set{2,4}\) have no upper bound.

  The algebra  \((\UNI,+,\oplus)\) of arithmetical functions with unitary
  convolution is defined as the set of all functions 
  \[\setN^+ \to \setC,\]
  with the structure of \(\setC\)-vector space given by
  component-wise addition and multiplication of scalars, and with
  multiplication defined by the 
  \emph{unitary convolution}
  \begin{equation}
    \label{eq:unitaryconv}
    ( f\oplus g) (n) =  \sum_{\ddivides{d}{n}} f(d) g(n/d) = \sum_{d
      \oplus m = n} f(d) g(m).
  \end{equation}
  \((\UNI,+,\oplus)\) is a non-Noetherian, quasi-local ring, where the
  non-units (i.e. the elements in the unique maximal ideal) consists
  of those \(f\) with \(f(1)=0\).

  It is
  natural to endow the algebra \(\UNI\) with the non-archimedean norm 
  \begin{equation}
    \label{eq:norm}
    \norm{f} = \frac{1}{\ord(f)}, \qquad \ord(f) = \min \supp(f),
  \end{equation}
where for \(f \in \UNI\),  \(\supp(f)=\setsuchas{n \in \setN^+}{f(n)
  \neq 0}\). If we give \(\setC\) the trivial norm 
\begin{displaymath}
  \norm{c} = 
  \begin{cases}
    1 & c \neq 0 \\
    0 & c= 0
  \end{cases}
\end{displaymath}
then \(\UNI\) becomes a normed vector space over \(C\), and
\((\UNI,+,\oplus)\) becomes a normed (not valued) algebra
\cite{Snellman:UniDivTop}. In detail: by a \emph{normed \(C\)-vector
  space} we mean a \(\setC\)-vector space \(A\) equipped with a norm, such that
for \(f \in A\), \(c \in \setC\) we have that
\begin{equation}
  \label{eq:normprop}
  \begin{split}
  \norm{f+g} & \leq \max (\norm{f},\norm{g})\\
  \norm{cf} &= \norm{c} \norm{f}
  \end{split}
\end{equation}
A \(\setC\)-algebra \(A\) is  normed if its underlying vector space is
normed, and if for \(f,g \in A\),
\begin{equation}
  \label{eq:normalg}
\norm{fg} \leq \norm{f} \norm{g}  
\end{equation}

  For \(n \in \setN^+\) we
  denote by \(e_n\) the function which is \(1\) on \(n\) and zero
  otherwise; we'll use the convention that \(e_0=0\) denotes the zero
  function. The function \(e_1\) is the multiplicative identity in
  \(\UNI\). 
  
  Every \(f\) in \(\UNI\) can be written uniquely as a
  convergent sum
  \begin{equation}
    \label{eq:Schauder}
    f = \sum_{n \in \setN^+} f(n) e_n,
  \end{equation}
  thus the \(e_n\)'s form a Schauder basis, and the set
  \[\setsuchas{e_{p^r}}{r \in \Nat, \, p \text{ prime }}\]
  is a minimal generating
  set in the sense that it generates a dense sub-algebra of \(\UNI\),
  and no proper subset has this property.

  \begin{lemma}\label{lemma:emult}
    For \(a,b \in \setN^+\), 
    \begin{displaymath}
      e_a \oplus e_b = e_{a \oplus b} =
      \begin{cases}
        e_{ab} & \text{ if } \gcd(a,b)=1 \\
        0 & \text{ otherwise } 
      \end{cases}
    \end{displaymath}
  \end{lemma}
  \end{subsection}

\begin{subsection}{The ring of arithmetical functions as a power series
    ring on a monoid-with-zero}
  
  \begin{definition}\label{def:monoidwzero}
    A monoid-with-zero \((M,0,1,\cdot)\) is a set \(M\), a zero element
    \(0 \not \in M\), a multiplicative unit \(1 \in M\), and a
    commutative, associative operation 
    \begin{displaymath}
       (M \cup \set{0}) × (M \cup \set{0}) \to M \cup \set{0}
    \end{displaymath}
     such that for all  \(m \in  M\),
     \begin{equation}\label{eq:monz}
       \begin{split}
         0 \cdot m &= 0 \\
         1 \cdot m &= m 
       \end{split}
     \end{equation}
 It is said to be cancellative as a monoid-with-zero if
    \begin{equation}\label{eq:canc}
      \forall a,b,c \in M: \,\, a \cdot b = a \cdot c \neq 0 \implies b = c.
    \end{equation}

A homomorphism 
\begin{equation}
  \label{eq:monzhom}
  f: (M,0,1,\cdot) \to (M',0',1',\cdot')
\end{equation}
of monoids-with-zero is a mapping \(M \cup \set{0} \to M' \cup
\set{0'}\) such that \(f(0)=0'\), \(f(1)=1'\) and \(f(x\cdot y) = f(x)
\cdot' f(y)\).

A subset \(S \subseteq M\) is a sub-monoid-with-zero if 
\[x,y \in S \implies x\cdot y \in S \cup \set{0},\]
 and a monoid ideal-with-zero if 
 \begin{displaymath}
   x \in S,  \, y \in M \quad \implies \quad x \cdot y \in S \cup \set{0}.
 \end{displaymath}
  \end{definition}

  \begin{remark}
    Dropping the demand for a unit \(1\), we get a
    \emph{semigroup-with-zero}. 
  \end{remark}

  \begin{lemma}
    \((\setN^+,\oplus)\) is a monoid-with-zero which is cancellative as
    a monoid-with-zero.
  \end{lemma}
  \begin{proof}
    We let define \(n \oplus 0 = 0\) for all \(n \in \setN^+\). Then
    \eqref{eq:monz} hold, with \(1 \in \setN^+\) as the multiplicative
    unit. If \(\gcd(a,b)=\gcd(a,c)=1\) and \(ab=ac\) then \(b=c\),
    so \((\setN^+,\oplus)\) is cancellative.
  \end{proof}

  We recall some of the definitions made in \cite{Snellman:UniDivTop}.
  Let \[Y=\setsuchas{y_{i,j}}{1 \leq i,j < \infty}\]
 be an \(\setN^+× \setN^+\)-indexed set of indeterminates, let
 \(Y^*\) denote the 
free abelian monoid on \(Y\), and let 
\begin{equation}
  \label{eq:separated}
  Y^* \supseteq \M = \set{1}  \cup \setsuchas{y_{i_1,j_1} \cdots
    y_{i_r,j_r}}{i_1 < \cdots < i_r}
\end{equation}
 denote the subset of \emph{separated monomials}. We regard \(\M\) as
 a monoid-with-zero, the multiplication given by 
 \begin{equation}
   \label{eq:sepmult}
   a \oplus b = 
   \begin{cases}
     ab & \text{ if } ab \in \M \\
     0 & \text{ otherwise}
   \end{cases}
 \end{equation}
Note that \((\M,0,1,\oplus)\) is not a sub-monoid-with-zero of
\((Y^*,0,1,\cdot))\), but it is an epimorphic image, under the map 
\begin{displaymath}
  Y \ni m \mapsto 
  \begin{cases}
    m & \text{ if } m \in \M \\
    0 & \text{ if } m \not \in \M
  \end{cases}
\end{displaymath}

Denote by \(\sqbij\)  the map between \(\M\) and
\(\setN^+\) 
determined by
\begin{equation}
  \label{eq:sqbij}
   \sqbij\left(\prod_{i=1}^r y_{a_i,b_i}\right) = \prod_{i=1}^r
   p_{a_i}^{b_i} 
\end{equation}

\begin{lemma}
  \(\sqbij\) induces an isomorphism of monoids-with-zero. 
\end{lemma}
\begin{proof}
  We extend \(\sqbij\) by putting \(\sqbij(0)=0\). If \(a,b \in \M\)
  then \(ab\) is separated if and only if \(\sqbij(a)\) and
  \(\sqbij(b)\) are relatively prime. Hence, \(\sqbij\) is a
  homomorphism. It follows from the fundamental theorem of arithmetic
  (unique factorization of positive integers) that \(\sqbij\) is a bijection.
\end{proof}

The following definition is similar to the more general
construction of Ribenboim \cite{Rib:GP,Rib:UniZ}.
\begin{definition}
  Let \((M,\cdot)\) be a commutative monoid (or let \((M,0,1,\cdot)\)
  be a monoid-with-zero) such that for any \(m \in M\), the
  equation 
  \begin{displaymath}
    x \cdot y = m
  \end{displaymath}
  have only finitely many solutions \((x,y) \in M × M\). Then the
  generalized power series ring on \(M\) with coefficients in
  \(\setC\), denoted \(\setC[[M]]\), 
  is the set of all maps \(f: M \to \setC\), with the obvious structure
  as a \(C\)-vector space, and with multiplication given by the
  convolution
  \begin{displaymath}
    f * g (m) = \sum_{x \cdot y = m} f(x) g(y).
  \end{displaymath}
  
  We give \(\setC[[M]]\) the topology of point-wise convergence, where
  \(\C\) is discretely topologized, so if \((f_v)_{v=1}^\infty\) is a
  sequence in  \(\setC[[M]]\) then
  \begin{equation}
    \label{eq:pointwise}
    f_v \to f \in \setC[[M]] \, \iff \, \forall m \in M: \exists
    V(m): \forall v > V(m): \, f_v(m) = f(m) 
  \end{equation}
\end{definition}

As an example, \(\setC[[Y^*]]\) is the ``large power series ring'' in the
sense of Bourbaki \cite{Bourbaki:Commutative} on the set of
indeterminates \(Y\).

The following easy result was proved in \cite{Snellman:UniDivTop}.
\begin{theorem}
  Let 
  \begin{equation}
    \label{eq:sqa}
    \SQA = \setsuchas{f \in \setC[[Y^*]]}{\supp(f) \cap \M = \emptyset}
  \end{equation}
  Then
  \(\SQA\) is a closed ideal in \(\setC[[Y^*]]\), and 
  the ideal  minimally generated by the
  set
  \begin{equation}
    \label{eq:mingensSQA}
    \setsuchas{y_{i,a} y_{i,b}}{i,a,b \in \setN^+, \, a \leq b}
  \end{equation}
has \(\SQA\) as its closure.
Furthermore,
  \begin{equation}
    \label{eq:sim}
    \UNI \simeq \setC[[(\setN^+,\oplus)]] \simeq \setC[[\M]] \simeq
    \frac{\setC[[Y^*]]}{\SQA} 
  \end{equation}
\end{theorem}
The isomorphisms in \eqref{eq:sim} are homeomorphisms, so the
norm-topology on \(\UNI\) coincides with the topology of point-wise
convergence.

\begin{definition}
  We let \(\setC[M] \subseteq \setC[[M]]\) denote the dense 
  sub-algebra of finitely supported maps \(M \to \setC\).
\end{definition}

We denote by \(\UNI^f\) the sub-algebra of finitely supported maps
\(\setN^+ \to \setC\). It follows that 
  \begin{equation}
    \label{eq:sim2}
    \UNI^f \simeq \setC[(\setN^+,\oplus)] \simeq \setC[\M] \simeq
    \frac{\setC[Y^*]}{\left( \setsuchas{y_{i,a} y_{i,b}}{i,a,b
          \in \setN^+, \, a \leq b} \right)}
  \end{equation}

\end{subsection}
\end{section}

\begin{section}{The truncations $\UNI_V$ --- definition and basic
    properties} 
  In what follows, \(V\) will, unless otherwise stated, denote a
  subset of \(\setN^+\).

  \begin{definition}
    Define
    \begin{equation}
      \label{eq:GammaV}
      \begin{split}
        \M_V &=\sqbij^{-1}(V) \\ 
      \UNI_V &= \setsuchas{f \in \UNI}{\supp(f) \subseteq V} \\
      \UNI_V^f &= \setsuchas{f \in \UNI^f}{\supp(f) \subseteq V}
      \end{split}
    \end{equation}
    With component-wise addition, and the modified multiplication
    \begin{equation}
      \label{eq:truncmult}
      (f \oplus_V g) (k) = 
      \begin{cases}
        (f \oplus g)(k) & k \in V \\
        0 & k \not \in V
      \end{cases}
    \end{equation}
    \(\UNI_V\) becomes a commutative algebra over \(\setC\). 
    The modified multiplication 
    \begin{equation}
      \label{eq:multV}
      a \oplus_V b = 
      \begin{cases}
        a\oplus b & a \oplus b \in V\\
        0 & a \oplus b \not \in V
      \end{cases}
    \end{equation}
    makes \((V,0,\oplus_V)\) into a semigroup-with-zero. We regard
    \(\M_V\) as a semigroup-with-zero
    isomorphic to 
    \((V,0,\oplus_V)\).

In what follows, we will often (by abuse of notation) use \(\oplus\)
for \(\oplus_V\). 

We define 
\begin{equation}
  \label{eq:smallerV}
  a \leq_{\oplus_V} c \quad \iff \quad \exists c: \, a \oplus_V b = c
\end{equation}

    \(\UNI_V\) and \(\UNI_V^f\) are unital, and \(\M_V\) is a
    monoid-with-zero, if and 
    only if \(1 \in V\).
  \end{definition}

  \begin{lemma}\label{lemma:closed}
    \(\UNI_V\) is a closed subspace and a sub vector space of
    \(\UNI\). If \(V\) is finite, then \(\UNI_V\) has the discrete
    topology, and \(\UNI_V=\UNI_V^f\). If \(V\) is infinite, then
    \(\UNI_V^f\) is dense in \(\UNI_V\).
  \end{lemma}
  \begin{proof}
    Suppose that \(f_v \to f\), \(f_v \in \UNI_V\), \(f \in \UNI
    \setminus \UNI_V\). Let \[j = \min\left( \supp(f) \setminus
      V)\right).\] Then \(\norm{f_v - f} \geq 1/j\), a contradiction. 
  \end{proof}

  \begin{lemma}
    \begin{equation}
      \label{eq:Vobv}
      \begin{split}
      \UNI_V &\simeq \setC[[(V,\oplus_V)]] \simeq \setC[[\M_V]]\\
      \UNI_V^f &\simeq \setC[(V,\oplus_V)] \simeq \setC[\M_V]
      \end{split}
    \end{equation}
    where the isomorphisms are also homeomorphisms.
  \end{lemma}

  \begin{definition}
Let \(\upsets\) denote the set of non-trivial sub-monoids of the
monoid-with-zero  \((\setN^+,\oplus)\). 
Thus \(W \in \upsets\) if and only if \(1\in W\) and  \(a,b \in W\) implies that \(a
\oplus b \in W \cup \set{0}\). 
  \end{definition}

  \begin{lemma}\label{lemma:whensub}
    Let \(1 \in W \subseteq \setN^+\). Then
    \(\UNI_W\) is a  closed sub-algebra of \(\UNI\) if and only if \(W
    \in \upsets\). 
  \end{lemma}
  \begin{proof}
    From Lemma~\ref{lemma:closed} we have that \(\UNI_W\) is a closed
    subspace of \(\UNI\), and a vector subspace. 
    So \(\UNI_W\) is a closed sub-algebra if and only if
    \((\UNI_W,\oplus)\) is a sub-monoid of the 
    multiplicative monoid-with-zero of \(\UNI\).

Suppose first that
    \(W \in \upsets\), and that 
    \(f,g  \in \UNI_W\). 
    To prove that \[\supp(f \oplus g) \subseteq W,\] take
    \(n \not \in W\). Since \(W\) is a sub-monoid of
    \((\setN,\oplus))\), \(n\) can not be written as \(n=a \oplus b\)
    with \(a,b \in W\). Thus 
    \begin{displaymath}
\left(f \oplus g\right)(n) = \sum_{a  \oplus b = n} f(a) g(b) = 0.      
    \end{displaymath}

    On the other hand, if \(W \not \in \upsets\), then there are \(a,b
    \in W\) with \(\gcd(a,b)=1\) such that \(a \oplus b \not \in W\). It
    follows that \[e_a \oplus e_b = e_{a \oplus b} \not \in \UNI_W.\]
  \end{proof}

  \begin{definition}
    For \(1 \in V \subseteq W \subseteq \setN^+\) we define the
    \emph{truncation map} 
    \begin{equation}
      \label{eq:rest}
      \begin{split}
      \rest{W,V}: \UNI_W & \to \UNI_V \\
      \rest{W,V}(f)(k) & = 
      \begin{cases}
        f(k) & k \in V \\
        0 & k \not \in V
      \end{cases}
      \end{split}
    \end{equation}
    We also define 
  \begin{equation}
    \label{eq:proj}
    \begin{split}
      \rest{V}: \UNI &\to   \UNI_V \\
    \rest{V}(f)(k) & = 
    \begin{cases}
      f(k) & k \in V \\
      0 & k \neq  V
    \end{cases}
    \end{split}
  \end{equation}
  \end{definition}

  \begin{lemma}
    The truncation maps \(\rest{V}\) and \(\rest{W,V}\) are
    surjective, \(\setC\)-linear 
    and continuous.  We have that
    \begin{equation}
    \label{eq:In}
    \ker \rest{V} = \GTRUNC_V = \setsuchas{f \in \UNI}{f(k) = 0 \text{
        for all } k \in V} = \UNI_{V^c}
  \end{equation}
  \(\GTRUNC_V\) is a closed subset of \(\UNI\).
  \end{lemma}
  \begin{proof}
    It is obvious that the maps are surjective and \(C\)-linear. Let us prove
    continuity for \(\rest{W,V}\). If \(f \in \UNI_W \setminus
    \set{0}\), then  
    \begin{displaymath}
      \begin{split}
      \norm{\rest{W,V}(f)} &= 
      \begin{cases}
        0 & \supp(f) \cap V = \emptyset \\
        \frac{1}{\min \setsuchas{k \in V}{f(k) \neq 0}}  & \text{ otherwise}
      \end{cases}
      \\
      & \leq \norm{f} = \frac{1}{\min \setsuchas{k \in W}{f(k) \neq 0}}
      \end{split}
    \end{displaymath}
    It follows that \(\rest{W,V}\) is continuous. The case of
    \(\rest{V}\) is 
    similar. Since \(\UNI\) is
    Hausdorff, so is \(\UNI_W\), and thus \(\ker \rest{V} =
    \rest{V}^{-1}(\set{0})\) is the 
    inverse image (under a continuous map) of a closed set, and
    consequently closed. 
  \end{proof}

  \begin{definition}
    Let \(\oksets\) denote the collection of non-empty order ideals of 
    \((\setN^+,\geq_{\oplus})\), and let \(\oksets^f\) denote the
    collection of finite, non-empty order ideals of 
    \((\setN^+,\geq_{\oplus})\). Thus \(V \in \oksets\) if and only if
    \begin{equation}
      \label{eq:oksetdef}
      n \in V, \, d \in \Nat^+, \, \ddivides{d}{n} \, \implies \, d
      \in V.
    \end{equation}
  \end{definition}

  \begin{lemma}\label{lemma:emptyorinfinite}
    If \(V \in \oksets\) then \(1 \in V\), and \(V^c\) is either empty
    or infinite. 
  \end{lemma}
  \begin{proof}
    If \(V=\setN^+\) then both assertions hold.
    If \(V \neq \setN^+\) then there is some \(n \in \setN^+ \setminus
    V\). Hence \(1\), which is a unitary divisor of \(n\), is in
    \(V\). Since \(V\) is an order ideal, it follows that whenever
    \(\gcd(n,m)=1\), then \(nm = n \oplus m \not \in V\). Thus the
    complement of \(V\) is infinite.
  \end{proof}

  \begin{theorem}
    Let \(1 \in V \subseteq \setN^+\).
    The following are equivalent:
    \begin{enumerate}[(i)]
    \item \(V \in \oksets\),
    \item \(V^c\) is a monoid ideal,
    \item \(\GTRUNC_V\) is a closed ideal of \(\UNI\), and
      \(\UNI_V = \frac{\UNI}{\GTRUNC_V}\).
    \item The truncation map \(\rest{V}\) is an algebra homomorphism.
    \end{enumerate}
  \end{theorem}
  \begin{proof}
    Suppose that \(V \in \oksets\). If \(x \in V^c\), \(y \in \setN^+\)
    then either \(x \oplus y = 0\) or \(x \oplus y \in V^c\), because
    otherwise \(x \oplus y \in V \implies x \in V\), a
    contradiction. So \(V^c\) is a monoid ideal.
    
    Let \(f,g \in   \UNI\),  \(n \in V\). Then 
  \begin{displaymath}
      \rest{V}(f \oplus g) (n) = 
        f \oplus h (n)  = \sum_{\ddivides{d}{n}} f(d) h(n/d)
    \end{displaymath}
    whereas 
    \begin{displaymath}
      \bigl[ \rest{V}(f) \oplus_V \rest{V}(h) \bigr] (n) =
        \sum_{\substack{\ddivides{d}{n}\\d \in V\\n/d \in V}} f(d)
        h(n/d)  = \sum_{\ddivides{d}{n}} f(d) h(n/d) 
    \end{displaymath}
    If \(m \not \in V\) then 
\begin{displaymath}
      \rest{V}(f \oplus g) (m) = \left[ \rest{V}(f) \oplus_V
        \rest{V}(g) \right] (m) = 0,
    \end{displaymath}
    so we have shown that 
    \(\rest{V}\) is a homomorphism, hence the kernel \(\GTRUNC_V\) is
    an ideal,     and \(\UNI_V = \frac{\UNI}{\GTRUNC_V}\).

    If \(V \not \in \oksets\) then there exist \(n \in V\), \(a,b \in
    \setN^+\), \(a \not \in V\), such that \(a \oplus b = n\). 
    Thus \(V^c\) is not a monoid ideal, and
    \begin{displaymath}
      \rest{V}(e_a \oplus e_b) = \rest{V}(e_n) = e_n \neq
      0 = \rest{V}(e_a) \oplus \rest{V}(e_b).
    \end{displaymath}
    So \(\rest{V}\) is not an algebra homomorphism. Furthermore, \(e_a
        \in \GTRUNC_V\), \(e_b \in \UNI\), but 
    \[e_a   \oplus e_b = e_n  \not \in \GTRUNC_V,\]
    so \(\GTRUNC_V\) is not an ideal.
  \end{proof}

  \begin{corr}
    Let \(1 \in V \subseteq \setN^+\).
    The following are equivalent:
    \begin{enumerate}[(i)]
    \item \label{en:retract}
      \(\UNI_V\) is a continuous algebra retract of \(\UNI\).
    \item \label{en:okup}
      \(V \in \oksets \cap \upsets\).
    \item \label{en:submo}
      \(V\) is a sub-monoid of \((\setN^+,\oplus)\) generated by a
      some subset of \(\PP\).
\end{enumerate}
  \end{corr}
  \begin{proof}
    Since a continuous homomorphism from \(\UNI\) to \(\UNI_V\) is
    determined by its values on the Schauder basis \(\setsuchas{e_n}{n
      \in \setN^+}\), \(\UNI_V\) is an algebra retract of \(\UNI\) if an
    only if it is a 
    sub-algebra and the restriction map \(\rest{V}: \UNI \to \UNI_V\)
    is an algebra  homomorphism. This yields the equivalence between
    (\ref{en:retract}) and (\ref{en:okup}).

    If \(V\) is a non-empty
    order ideal in 
    \((\setN^+, \geq_\oplus)\), then it contains along with any \(n \in
    V\) all of its prime power unitary divisors, hence if it is also a
    monoid ideal in \((\setN^+,\oplus)\), it must contain the sub-monoid
    generated by those prime powers. Conversely, a  sub-monoid of
    \((\setN^+,\oplus)\) generated by a some subset of \(\PP\)
    is an order ideal. So (\ref{en:okup}) and (\ref{en:submo}) are
    equivalent. 
  \end{proof}

  \begin{definition}\label{definition:squarefree}
    Let \(\sqfree\) denote the set of positive, square-free integers;
    in particular, \(1 \in \sqfree\).
  \end{definition}

  \begin{lemma}\label{lemma:isOK}
    Let \(n \in \setN^+\). Then the following hold:
    \begin{enumerate}[(i)]
    \item \label{en:n1} \([n] \in \oksets^f\).
    \item \label{en:st1} \(\setsuchas{q \in \setN^+}{\ddivides{q}{n}} \in
      \oksets^f\).
    \item \label{en:sqf} \(\sqfree \in \oksets \cap \upsets\).
    \item \label{en:n2} \([n] \cap \sqfree \in \oksets^f\).
    \item \label{en:st2} 
      \(\setsuchas{q \in \sqfree}{\ddivides{q}{n}}   \in  \oksets^f\).
    \end{enumerate}
  \end{lemma}
  \begin{proof}
    \eqref{en:n1}: If \(\ddivides{w}{n}\) then \(w \le n\) so \(w \in
    [n]\).
    
    \eqref{en:st1}: In any locally finite poset with a minimal
    element, the principal order ideal on an element is a finite order 
    ideal.

    \eqref{en:sqf}: The  divisors of a square-free integer are
    square-free, so in particular the unitary divisors are
    square-free. If \(a,b\) are square-free then \(a \oplus b\) is
    either square-free (if \(\gcd(a,b)=1\)) or zero.

    \eqref{en:n2}, \eqref{en:st2}: Intersections of order ideals are
    order ideals.
  \end{proof}

  \begin{lemma}
    Let \(1 \in V \subseteq W \subseteq \setN^+\), \(V \in \oksets\). Then
    \(\rest{W,V}\) is an 
    algebra epimorphism.
  \end{lemma}
  \begin{proof}
    Let \(f,g \in \UNI_W\). If \(n \in W \setminus V\) then 
    \begin{displaymath}
      0=\rest{W,V}(f \oplus_W g)(n) = \bigl(\rest{W,V}(f) \oplus_V
      \rest{W,V}(g)\bigr )(n) =0,
    \end{displaymath}
    whereas if \(n \in V\)  then
    \begin{displaymath}
      \begin{split}
      \rest{W,V}(f \oplus_W g)(n) &= (f \oplus_W g)(n) \\
      &=  \sum_{\substack{\ddivides{d}{n}\\ d \in W}} f(d) g(n/d) \\
      &=  \sum_{\substack{\ddivides{d}{n}\\ d \in V}} f(d) g(n/d) \\
      \end{split}
    \end{displaymath}
    where the last equality follows since 
     \(V\) is an order ideal in \(W\), and consequently  any unitary divisor of
    \(n\) is in \(V\), so \(d, n/d \in V\).
Using this observation, we have
    \begin{align*}
      \bigl(
      \rest{W,V}(f) \oplus_V  \rest{W,V}(g)
      \bigr)(n) &= \sum_{\substack{\ddivides{d}{n}\\ d \in V}}
      \rest{W,V}(f) (d) \rest{W,V}(g) (n/d) \\
      &= \sum_{\substack{\ddivides{d}{n}\\ d \in V}}
      f(d) g (n/d) \\
      &= \rest{W,V}(f \oplus_W g)(n)
    \end{align*}
  \end{proof}

  \begin{corr}\label{corr:cyclicUNIW}
    If \(1 \in V \subseteq W \subseteq \setN^+\), \(V \in \oksets\), then
    \(\UNI_V\) is a cyclic \(\UNI_W\)-module.
  \end{corr}

The following theorem gives a motivation for our studies of the
truncations \(\UNI_V\): they approximate \(\UNI\) in the natural
sense.

  \begin{theorem}\label{thm:approx}
    Let \(\finsubs\) denote the set of finite subsets of \(\Nat^+\).
    \begin{enumerate}[(A)]
    \item \label{en:allV}
    The set of all \(\UNI_V\),  \(V \in \finsubs\), and truncation maps 
    \[\rest{W,V}: \UNI_W  \to \UNI_V, \qquad  V ,W \in \finsubs, \, V
    \subseteq W
    \]
    forms an inverse system
    of normed vector spaces over \(\setC\), and
    \begin{equation}
      \label{eq:invsys}
    \varprojlim_{V \in \finsubs} \UNI_V \simeq \UNI      
    \end{equation}
    as normed vector spaces over \(\setC\).

  \item \label{en:okcofinal}
\(\oksets^f\) is cofinal in  \(\finsubs\), and
    \begin{equation}
      \label{eq:invsys2}
    \varprojlim_{V \in \oksets^f} \UNI_V \simeq \UNI      
    \end{equation}
    as normed vector spaces over \(\setC\).

  \item \label{en:okalg}
    The  set of all \(\UNI_V\), \(V \in \oksets^f\),
    together with  all truncation 
    maps 
    \[\rest{W,V}: \UNI_W \to \UNI_V, \qquad V \subseteq W, \, V,W \in
    \oksets^f,\] 
    forms an inverse system of normed Artinian \(\setC\)-algebras, and
    \eqref{eq:invsys2} is an isomorphism of normed \(\setC\)-algebras.
    \end{enumerate}
  \end{theorem}

  \begin{proof}
    \eqref{en:allV}
    Let \(Q\) be a normed \(\setC\)-vector space, and consider
    the diagram 
    \begin{equation}
      \label{eq:invlim}
      \xymatrix{
        \UNI \ar [rd]_{\rest{V}} \ar @/_1pc/ [rdd]_{\rest{W}}& & Q
            \ar @{.>} [ll]_{g}
            \ar [ld]^{f_V}
            \ar @/^1pc/ [ldd]^{f_W}
            \\ 
        & \UNI_V & \\
        & \UNI_W \ar [u]_{\rest{W,V}}&
        }
    \end{equation}
    where \(W \supseteq V\) and \(f_W\), \(f_V\) are given continuous
    homomorphisms. If the diagram without the dotted line
    commutes, the whole diagram commutes when we define 
    \begin{equation}
      \label{eq:g}
      g(x)(k) = f_{V'}(x)(k), \qquad V' \ni k
    \end{equation}
    Clearly, \(g\) is \(\setC\)-linear; we claim that it  is  also
    continuous.
    To prove this, let \(x_n \to 0\) in \(Q\).
    Then \(f_V(x_n) \to 0\) in every \(\UNI_V\), that is, 
    \(\norm{f_V(x_n)} \to 0\). Fix a \(k \in \setN^+\), and suppose
    that \(k \in V\). Since \(\norm{f_V(x_n)} \to 0\), there is an
    \(N\) such that \(\norm{f_V(x_n)} < 1/k\) whenever \(n >
    N\). Hence for \(n > N\)
    \begin{displaymath}
      \norm{g(x_n)} = \norm{f_V(x_n)} < 1/k.
    \end{displaymath}

    \eqref{en:okcofinal}: If \(V \subseteq \Nat^+\) is finite, so is the
    set obtained by adding all unitary divisors of elements in
    \(V\). This proves that \(\oksets^f\) is cofinal in \(\finsubs\).

    \eqref{en:okalg}
    Let \(Q\) be a normed algebra, and consider
    once again the  diagram \eqref{eq:invlim}, where now the \(f_V\)'s
    are continuous algebra 
    homomorphisms.
    We must show that now \(g\) is also a ring
    homomorphism. So let \(x,y \in Q\), fix \(k \in \setN^+\), 
    and take \(V \in \oksets^f\)    with \(k \in V\). Then
    \begin{displaymath}
      \begin{split}
      g(xy)(k) &= f_V(xy)(k) \\
      & = \bigl( f_V(x) \oplus_V f_V(y) \bigr) (k) \\
      & = \sum_{\ddivides{d}{k}} f_V(x)(d) f_V(y)(k/d) \\ 
      & = \sum_{\ddivides{d}{k}} g(x)(d) g(y)(k/d) \\ 
      &= \bigl( g(x) \oplus_V g(y) \bigr) (k),
      \end{split}
    \end{displaymath}
    where we have used that \(f_V\) is a ring homomorphism, and that
    \(k \in V \in \oksets^f\) implies that any unitary divisor of \(k\)
    lies in \(V\).

    Hence, \(g(xy) = g(x) \oplus g(y)\).
  \end{proof}

With minor modifications, the preceding proof works for the following
generalization: 
  \begin{theorem}
    Let \(U \subseteq \setN^+\).
    \begin{equation}
      \label{eq:minvsys}
    \varprojlim_{\substack{V \subseteq U\\ V \in \finsubs}} \UNI_V \simeq \UNI_U      
    \end{equation}
    as normed vector spaces over \(\setC\), and
    \begin{equation}
      \label{eq:minvsys2}
    \varprojlim_{\substack{V \in \oksets^f\\V \subseteq U}}
    \UNI_V \simeq  \UNI_U
    \end{equation}
    as normed algebras over \(\setC\).
  \end{theorem}

    \begin{lemma}\label{lemma:units}
      Let \(V \in \oksets\).
      \begin{enumerate}[(i)]
      \item       The multiplicative unit of \(\UNI_V\) is
        \(\rest{V}(e_1)\). 
      \item \label{en:ql}
        \(\UNI_V\) is quasi-local: \(f \in \UNI_V\) is a unit if and only if
      \(f(1)\neq 0\), and the non-units form the unique maximal ideal.
    \item       If \(V\) is finite, then all non-units are nilpotent.
    \item \(\UNI_V\) is Artinian if and only if \(V\) is finite.
    \item \label{en:noe}
      \(\UNI_V\) is Noetherian if and only if \(V\) is finite.
      \end{enumerate}
    \end{lemma}
    \begin{proof}
      We'll prove (\ref{en:ql}) and (\ref{en:noe}),
      the rest is trivial.

      (\ref{en:ql}): 
      If \(f(1) \neq 0\), then let
      \(f'\) denote any ``lift'' of \(f\) to \(\UNI\), i.e.
      \(\rest{V}(f')=f\). Then there is a \(g' \in \UNI\) such that
      \(f'g' = e_1\), hence 
      \begin{displaymath}
        \rest{V}(f'g') = f \rest{V}(g') =\rest{V}(e_1).
      \end{displaymath}
      Thus \(f\) is a unit, hence regular.

      Conversely, it is easy to see that if \(f(0)=0\) then \(f\) is
      not a unit. If \(V\) is finite, let \(n = \max V\), and let \(g
      \in \UNI\) be such 
      that \(\rest{V}(g)=f\). We claim
      that if \(a,b \in \UNI\) are (non-zero) non-units, then 
      \[\ord(a \oplus b) \geq \ord(a) \ord(b) > \max
      \set{\ord(a),\ord(b)}.\]  Thus we 
      conclude that
       \(\ord(g^{n+1}) \geq n+1\), 
      which implies that \[\rest{n}(g^{n+1}) = f^{n+1}= 0.\]

      To prove the claim, let \(i=\ord(a)\), \(j=\ord(b)\). Then for
      \(k < ij\), 
      \[a\oplus b(k) = \sum_{\ddivides{\ell}{k}} a(\ell)
        b(k/\ell) = 0,\] since it is impossible that \(\ell \geq i\)
        and \(k/\ell \geq j\).

        (\ref{en:noe}):  If \(V\) is finite, then \(\UNI_V\) is
        Artinian, hence Noetherian. If \(V\) is infinite, then since
        \(V \in \oksets\) we must have that \(V \cap \PP\) is
        infinite, as well. The ideal generated by 
        \begin{displaymath}
          \setsuchas{e_q}{q \in V \cap \PP}
        \end{displaymath}
        is not finitely generated.
    \end{proof}

  \begin{theorem}\label{thm:mingens}
    Let \(V \in \oksets\). The set 
   \begin{equation}
     \label{eq:mingensMV}
     M_V = \setsuchas{e_k}{k \not \in V, \text{ but } d \in V \text{
         for all proper unitary divisors } d \text{ of } k}
   \end{equation}
   form a minimal
   generating set of an ideal \(I_V\) whose closure is \(\GTRUNC_V\).
  \end{theorem}
  \begin{proof}
    Recall that \(\setsuchas{e_n}{n \in \setN^+}\) is a Schauder base
    for \(\UNI\). Thus if  
    \begin{displaymath}
      \UNI \ni f = \sum_{n=1}^\infty c_n e_n, \qquad f \in \GTRUNC_V
    \end{displaymath}
    then 
    \begin{equation}\label{eq:sss}
      f = \sum_{n \not \in V} c_n e_n
    \end{equation}
    It follows that 
    \begin{displaymath}
      \setsuchas{e_n}{n \not \in V}
    \end{displaymath}
    forms a generating set of an ideal \(I_V\) s.t. \(I_V
    \subseteq\GTRUNC_V = \overline{I_V}\). Clearly, \(M_V\) is a
    minimal generating set of \(I_V\).
  \end{proof}

  \begin{lemma}
    Let \(V \in \oksets\), \(V \neq \setN^+\). Then
     \(I_V\) is not
    finitely generated, and \(I_V \subsetneq \GTRUNC_V\).
  \end{lemma}
  \begin{proof}
From the proof of Lemma~\ref{lemma:emptyorinfinite} we have that 
there exists \(n   \in V^c\), \(q_1,q_2,q_3,\dots \in
\Primes\)  such that \((n,q_i)=1\) and \(q_i, nq_i \in V^c\).  We can
furthermore assume that all unitary divisors of \(n\) belong to
\(V\). Then \(e_{nq_i}\) is a minimal generator of \(I_V\). 
The sums of these \(e_{nq_i}\)'s is an element in \(\GTRUNC_V\)
  but not (since the sum is infinite) in \(I_V\).
  \end{proof}

  \begin{theorem}\label{thm:pres}
    Let \(V \in \oksets\). Then there is a smallest \(W \in \oksets
    \cap \upsets\) containing \(V\), and a subset \(P \subseteq
    \PP\) such that \(W\) is the sub-monoid of
    \((\setN^+,\oplus)\) generated by \(P\). 
    Define
  \begin{equation}
    \label{eq:Y(V)}
    \begin{split}
    Y \supseteq Y(V) &= \sqbij^{-1}(\bar{V}) \cap Y \\
    &= \setsuchas{y_{i,j}}{\exists v \in V: \, \ddivides{p_i^j}{v}}\\
    &= \setsuchas{y_{i,j}}{p_i^j \in V}
    \end{split}
  \end{equation}
    and denote (by abuse of notation) by \(\setC[[Y(V)]]\) the
    generalized power series ring on the free abelian sub-monoid of
    \(Y^*\) which \(Y(V)\) generates. Similarly, denote by 
    \(\setC[Y(V)]\) the polynomial ring on the set of variables \(Y(V)\).
    Then \(\UNI_V\) is
    a cyclic \(\setC[[Y(V)]]\)-module, as well as a cyclic
    \(\UNI_W\)-module; \(\UNI_V^f\) is
    a cyclic \(\setC[Y(V)]\)-module, as well as a cyclic
    \(\UNI_W^f\)-module; 
  \end{theorem}
  \begin{proof}
    Let \(P=V \cap  \PP\), and let \(W\) be the sub-monoid it
    generates. Then \(W \in \oksets\), proving the first assertion.
    By Corollary~\ref{corr:cyclicUNIW} we have that \(\UNI_V\) is a cyclic
    module over \(\UNI_W\), which in turn is a cyclic module over
    \(\UNI\), which in turn is a cyclic module over \(\setC[[Y^*]]\). So
    \(\UNI_V\) is a cyclic \(\setC[[Y^*]]\)-module. Since every
    variable 
    not in \(Y(P)\) will act trivially on \(\UNI_V\), the latter is in
    fact a cyclic \(\setC[[Y(P)]]\)-module.
  \end{proof}

\end{section}

\begin{section}{Hilbert and Poincaré-Betti series}
    \begin{definition}
    We denote by \(\setN^\omega\) the abelian monoid of all finitely
    supported functions
    \(\setN^+ \to \setN,\) with component-wise addition. A
    \(\setC\)-vector 
    space \(R\) is \(\omega\)-multi-graded it is graded over
    \(\setN^\omega\). The grading is \emph{locally finite} if
    \(R_\vektor{\alpha}\) is a finite dimensional for all
    \(\vektor{\alpha} \in \setN^\omega\).
    For such an \(R\), we define the \(\setN^\omega\)-graded Hilbert
    series of \(R\) as 
    \begin{equation}
      \label{eq:hilbert}
      R(\vektor{u}) = \sum_{\vektor{\alpha} \in \setN^\omega}
      \dim_\setC 
      R_{\vektor{\alpha}} \vektor{u}^\vektor{\alpha} \in
      \Z[[\vektor{u}]] 
    \end{equation}
    If \(c: \setN^\omega \to \setN^s\)
    is a monoid homomorphism, then an \(\omega\)-multi-graded \(R\) is 
    \(s\)-multi-graded by
    \begin{equation}
      \label{eq:omegagraded}
      R_{\vektor{\gamma}} := \bigoplus_{\substack{\vektor{\alpha} \in
        \setN^\omega\\ c(\vektor{\alpha}) = \vektor{\gamma}}}
    R_\vektor{\alpha} 
    \end{equation}
    We say that the \(s\)-multi-grading is obtained by
    \emph{collapsing} the \(\setN^\omega\)-grading. 
    \(c\) is called locally finite if \(c^{-1}(\vektor{\gamma})\) is
    finite for all \(\vektor{\gamma} \in \setN^s\). Hence, if the
    \(\setN^\omega\)-grading is locally finite then the 
    collapsed grading is locally finite if and only if \(c\) is locally finite.
    Then, we define
    \begin{equation}
      \label{eq:hilbert2}
      \begin{split}
      R(t_1,\dots,t_s) &= R(c(u_1),c(u_2),\dots)\\
      &= \sum_{\vektor{\gamma} \in \setN^s} \dim_\setC
      R_{\vektor{\gamma}} \vektor{t}^\vektor{\gamma} \in
      \Z[[t_1,\dots,t_s]]         
      \end{split}
    \end{equation}
  \end{definition}

  \begin{definition}
    Define a bijection
    \begin{equation}
      \label{eq:nu}
      \begin{split}
      \nu: Y & \to \setN^+ \\
    \nu(y_{i',j'}) > \nu(y_{i,j}) & \iff \quad
    \Phi(y_{i',j'}) >  \Phi(y_{i,j})
      \end{split}
    \end{equation}
    Then \(\nu\) extends to a monoid isomorphism \(Y^* \to
    \setN^\omega\), such that
    \[y_{i,j}^k \mapsto  (0,\dots,0,k,0,\dots) \in \setN^\omega,\]
    with \(k\) in the \(\nu(y_{i,j})\)'th position.
    Hence, we may regard a \(\setN^\omega\)-graded algebra as
    \(Y^*\)-graded. Furthermore, we'll call an \(Y^*\)-graded algebra,
    which is concentrated in the degrees \(\M \subseteq Y^*\),
    \(\M\)-graded; a \(\setN^\omega\)-graded algebra concentrated in
    the corresponding multi-degrees of \(\Nat^\omega\) will also be
    called \(\M\)-graded.

  \end{definition}

  The enumeration of \(Y\) given by \eqref{eq:nu} is illustrated in
  Table~\ref{tab:nuenum}, as well as the induced
  enumeration for the variables corresponding to square-free elements.

    \begin{table}[htbp]
    \setlength{\unitlength}{0.8cm}
    \begin{tabular}{|c|c|}
    \begin{minipage}[c]{6cm}
    \begin{picture}(6,7)(-1,-1)
      \multiput(0.3,0)(0,1){6}{\multiput(0,0)(1,0){6}{\circle*{0.1}}}
      \put(0.0,0.2){\(y_{1,1}\)}
      \put(0.0,1.2){\(y_{1,2}\)}
      \put(0.0,2.2){\(y_{1,3}\)}
      \put(0.0,3.2){\(y_{1,4}\)}
      \put(1.0,0.2){\(y_{2,1}\)}
      \put(1.0,1.2){\(y_{2,2}\)}
      \put(1.0,2.2){\(y_{2,3}\)}
      \put(1.0,3.2){\(y_{2,4}\)}
      \put(2.0,0.2){\(y_{3,1}\)}
      \put(2.0,1.2){\(y_{3,2}\)}
      \put(2.0,2.2){\(y_{3,3}\)}
      \put(2.0,3.2){\(y_{3,4}\)}
      \put(3.0,0.2){\(y_{4,1}\)}
      \put(3.0,1.2){\(y_{4,2}\)}
      \put(3.0,2.2){\(y_{4,3}\)}
      \put(3.0,3.2){\(y_{4,4}\)}

      \put(-0.2,2.8){\line(1,0){1}}
      \put(0.8,2.8){\line(0,-1){1}}
      \put(0.8,1.8){\line(1,0){1}}
      \put(1.8,1.8){\line(0,-1){1}}
      \put(1.8,0.8){\line(1,0){2}}
      \put(3.8,0.8){\line(0,-1){1.3}}
      \put(3.8,-0.5){\line(-1,0){4}}
      \put(-0.2,-0.5){\line(0,1){3.3}}
      
    \end{picture}

    \end{minipage}
    &

  \begin{minipage}[c]{6cm}
    \begin{picture}(6,7)(-1,-1)

      \multiput(0.3,0)(0,1){6}{\multiput(0,0)(1,0){6}{\circle*{0.1}}}

      \put(0.0,0.2){2}
      \put(0.0,1.2){4}
      \put(0.0,2.2){8}
      \put(0.0,3.2){16}
      \put(1.0,0.2){3}
      \put(1.0,1.2){9}
      \put(1.0,2.2){27}
      \put(1.0,3.2){81}
      \put(2.0,0.2){5}
      \put(2.0,1.2){25}
      \put(2.0,2.2){125}
      \put(2.0,3.2){625}
      \put(3.0,0.2){7}
      \put(3.0,1.2){49}
      \put(3.0,2.2){}
      \put(3.0,3.2){}
      \put(4.0,0.2){11}
      \put(5.0,0.2){13}

      \put(-0.2,2.8){\line(1,0){1}}
      \put(0.8,2.8){\line(0,-1){1}}
      \put(0.8,1.8){\line(1,0){1}}
      \put(1.8,1.8){\line(0,-1){1}}
      \put(1.8,0.8){\line(1,0){2}}
      \put(3.8,0.8){\line(0,-1){1.3}}
      \put(3.8,-0.5){\line(-1,0){4}}
      \put(-0.2,-0.5){\line(0,1){3.3}}

    \end{picture}

  \end{minipage}
\\ \hline

  \begin{minipage}[c]{6cm}
    \begin{picture}(6,7)(-1,-1)

      \multiput(0.3,0)(0,1){6}{\multiput(0,0)(1,0){6}{\circle*{0.1}}}

      \put(0.0,0.2){1}
      \put(0.0,1.2){3}
      \put(0.0,2.2){6}
      \put(0.0,3.2){10}
      \put(1.0,0.2){2}
      \put(1.0,1.2){7}
      \put(1.0,2.2){}
      \put(1.0,3.2){}
      \put(2.0,0.2){4}
      \put(2.0,1.2){}
      \put(2.0,2.2){}
      \put(2.0,3.2){}
      \put(3.0,0.2){5}
      \put(3.0,1.2){}
      \put(3.0,2.2){}
      \put(3.0,3.2){}
      \put(4.0,0.2){8}
      \put(5.0,0.2){9}

      \put(-0.2,2.8){\line(1,0){1}}
      \put(0.8,2.8){\line(0,-1){1}}
      \put(0.8,1.8){\line(1,0){1}}
      \put(1.8,1.8){\line(0,-1){1}}
      \put(1.8,0.8){\line(1,0){2}}
      \put(3.8,0.8){\line(0,-1){1.3}}
      \put(3.8,-0.5){\line(-1,0){4}}
      \put(-0.2,-0.5){\line(0,1){3.3}}

    \end{picture}

  \end{minipage}
  & 
      \\ \hline
    \begin{minipage}[c]{6cm}
    \begin{picture}(7,2)(-1,-1)
      \multiput(0,0)(1,0){7}{\circle*{0.1}}
      \put(0.0,0.2){\(y_{1,1}\)}
      \put(1.0,0.2){\(y_{2,1}\)}
      \put(2.0,0.2){\(y_{3,1}\)}
      \put(3.0,0.2){\(y_{4,1}\)}
    \end{picture}

    \end{minipage}
    &

  \begin{minipage}[c]{6cm}
    \begin{picture}(7,2)(-1,-1)
      \multiput(0,0)(1,0){7}{\circle*{0.1}}
      \put(0.0,0.2){2}
      \put(1.0,0.2){3}
      \put(2.0,0.2){5}
      \put(3.0,0.2){7}
      \put(4.0,0.2){11}
      \put(5.0,0.2){13}
    \end{picture}

  \end{minipage}

\\ \hline
  \begin{minipage}[c]{6cm}
    \begin{picture}(7,2)(-1,-1)

      \multiput(0,0)(1,0){7}{\circle*{0.1}}
      
      \put(0.0,0.2){1}
      \put(1.0,0.2){2}
      \put(2.0,0.2){4}
      \put(3.0,0.2){5}
      \put(4.0,0.2){8}
      \put(5.0,0.2){9}
    \end{picture}

  \end{minipage}
&
\end{tabular}
        \caption{\(Y\) and the corresponding prime powers, and their
          enumeration. \(Y([10])\) is 
          enclosed. The square-free part of \(Y\) and the corresponding
          prime powers, and their 
          enumeration.} 
        \label{tab:nuenum}
    \end{table}

    We see that \(\setC[Y^*]\) is \(Y^*\)-graded, hence
    \(\UNI^f=\frac{\setC[Y^*]}{(y_{i,j}y_{i,k})}\) is also
    \(Y^*\)-graded; in fact, it is \(\M\)-graded. We have that
      \begin{equation}
        \label{eq:chara}
        \begin{split}
          \setC[Y^*](\vektor{u}) & = \sum_{m \in Y^*} m = \prod_{i,j=1}^\infty
          \frac{1}{1-u_{i,j}} \in \Z[[u_{i,j}]]\\
          \UNI^f(\vektor{u}) &= \prod_{i=1}^\infty
          \left(1+ \sum_{j=1}^\infty u_{i,j}\right) \\
        &1 + u_{1,1} + u_{2,1} + u_{1,2} +
          u_{3,1} + u_{1,1} u_{2,1} + \cdots \in \Z[[u_{i,j} \, \lvert
          \, i,j \in \Nat^+]]\\
          \UNI_{\sqfree}^f(\vektor{u}) &= 
          \prod_{i=1}^\infty
          \left(1+  u_{i,1} \right) \\
          &= 1 + u_{1,1} + u_{2,1} + u_{3,1} + u_{1,1}u_{2,1} + \cdots
          \in
          \Z[[u_{1,j} \, \lvert j \in \Nat^+]]
        \end{split}
      \end{equation}
or, regarding instead these algebras as \(\setN^\omega\)-graded,
\begin{equation}
  \label{eq:chara2}
  \begin{split}
          \setC[Y^*](\vektor{t}) & = \sum_{\vektor{\alpha} \in
            \setN^\omega} \vektor{t}^\vektor{\alpha} =
          \prod_{i=1}^\infty     \frac{1}{1-t_i}  \in \Z[[t_1,t_2,t_3,\dots]]
          \\
        \UNI^f(\vektor{t}) 
        &=\sum_{i=1}^\infty   \vektor{t}^{\mathrm{multideg}(e_i)} \\
        &= 1 + t_1 + t_2 + t_3 +t_4 + t_1t_2 + \dots
       \in \Z[[t_1,t_2,t_3,t_4,t_5,\dots]]
       \\
       \UNI_{\sqfree}^f(\vektor{t}) 
        &=\sum_{i \in \sqfree}  \vektor{t}^{\mathrm{multideg}(e_i)} \\
        &= 1 + t_1 + t_2 +t_4 + t_1t_2 + \dots
       \in \Z[[t_1,t_2,t_4,t_5, t_8, t_9,\dots]]
      \end{split}
\end{equation}

Note that the collapsing \(\M  \to \setN\) obtained by giving each
\(y_{i,j}\) degree 1 is \emph{not} locally finite, so we can not
define the ordinary \(\setN\)-graded Hilbert series of \(\UNI^f\).
There are however other ways of collapsing the grading that work:

\begin{example}\label{example:pairwise}
  Define a locally finite \(\setN^2\)-grading on \(\setC[Y^*]\) by 
    giving \(y_{i,j}\) bi-grade \((1,p_i^j)\). Then
    \begin{equation}
      \label{eq:gold}
      \begin{split}
       \UNI^f(u_1,u_2) 
       &=\prod_{i=1}^\infty
          \left(1+ \sum_{j=1}^\infty u_1 u_2^{p_i^j}\right) \\
          & =1 +\sum_{i,j=1}^\infty d(i,j) u_1 u_2^j \\
          &= 1+ \sum_{j=1}^\infty u_2^j \left( \sum_{i=1} d(i,j) u_1^i
          \right) \\
          &= 1+ u_1u_2^2 + u_1u_2^3 + u_1u_2^4 + (u_1+u_1^2)u_2^5 + (u_1+2u_1^2)u_2^7 + \dots
      \end{split}
    \end{equation}
    where \(d(i,j)\) is the number of ways the positive integer
    \(j\) 
    can be written as a sum of \(i\) pairwise co-prime prime powers.
    Collapsing this grading further, so that \(y_{i,j}\) gets degree
    \(p_i^j\), we have that
      \begin{equation}\label{eq:goldc}
        \begin{split}
       \UNI^f(u) &= \prod_{i=1}^\infty
          \left(1+ \sum_{j=1}^\infty u^{p_i^j}\right) \\
          & = \sum_{i=1}^\infty g(i) u^i \\
          &= 1+ u^2 + u^3 + u^4 + 2u^5 + 3u^7 + \dots
        \end{split}
      \end{equation}
    where \(g(i)\) is the number of ways the positive integer
    \(i\)   can be written as a finite sum of  pairwise co-prime prime
    powers. 
\end{example}

\begin{example}
  If \(R\) denotes the ring of finitely supported arithmetical
  functions \emph{with Dirichlet convolution}, then 
  \begin{displaymath}
    R \simeq \C[x_1,x_2,x_3,\dots],
  \end{displaymath}
  by letting \(e_n\), for \(n=p_{i_1}^{a_1} \cdots p_{i_r}^{a_r}\),
  correspond to the monomial \(x_{i_1}^{a_1} \cdots x_{i_r}^{a_r}\).
  We have that \(R\) is \(\Nat^\omega\)-graded, with Hilbert series
  \begin{displaymath}
    \begin{split}
    R(\vektor{u}) &= \sum_{\vektor{\alpha} \in \Nat^\omega}
    \vektor{u}^{\vektor{\alpha}}\\
    &= \prod_{i=1}^\infty \frac{1}{1-u_i}
    \end{split}
  \end{displaymath}
  Giving \(e_n\) the same bi-grade as in the previous example means
  giving \(x_i\) bi-grade \((1,p_i)\). Then 
  \begin{displaymath}
    \begin{split}
      R(u_1,u_2) &= \prod_{i=1}^\infty \frac{1}{1-u_1u_2^{p_i}}\\
      & = \sum_{i,j} c_{i,j} u_1^i u_2^j
    \end{split}
  \end{displaymath}
  where 
  \(c_{i,j}\) is the number of ways of writing \(j\) as the sum of
  \(i\) primes. Thus, by Vinogradov's three-primes theorem,
  \(c_{i,3}>0\) for odd \(i 
  \gg 0\), and if Goldbach's conjecture is true, then \(c_{i,3}>0\) for 
  \(i > 5\), and \(c_{i,2}>0\) for  odd \(i > 4\).

  Comparing this to \eqref{eq:gold}, one can ask if there is an \(j\)
  such that \(d_{i,j}>0\) for all sufficiently large \(i\), and what
  the minimal such \(j\) might be. We have checked that for any integer
  \(6 < i < 485\), either \(i \in \PP\) and \(d_{i,1}>0\), or 
  \(d_{i,2}>0\), thus one might conjecture that any integer \(>6\) is
  either a prime power or can be written as the sum of two relatively
  prime prime powers.
\end{example}

  \begin{definition}
    For an \(\setN^r\)-graded \(\setC\)-algebra \(R\) and a
    \(\setN^r\)-graded \(R\)-module \(M\), we
    regard \(\setC\) as a cyclic \(R\)-module, and
    define the \(\setN^r\)-graded \emph{Betti numbers} as
\begin{equation}
  \label{eq:mbetti}
  \beta_{i,\vektor{a}}(R,M) = \dim_\setC
  \Tor_R^i(M,\setC)_{\vektor{a}}
\end{equation}
and the     \emph{Poincaré-Betti series} by
\begin{equation}
  \label{eq:poincare}
  P_R^M (t,\vektor{u}) = \sum_{i=0}^\infty \sum_{\vektor{a} \in
    \setN^r}  \beta_{i,\vektor{a}}(R,M) t^i {\vektor{u}^{\vektor{a}}} 
  \end{equation}
We use the convention \(P_R(t,\vektor{u}) = P_R^\setC(t,\vektor{u})\).
  \end{definition}

If \(R\) and \(M\) are instead locally finite and
\(\setN^\omega\)-graded, then each \(\Tor_R^i(M,\setC)\) will be locally
finite and 
\(\setN^\omega\)-graded, so we can define the \(\setN^\omega\)-graded, 
or \(Y^*\)-graded,  
Poincaré-Betti series. 
That means that 
\begin{displaymath}
  P_R^M (t,\vektor{u}) \in \Z[[t]][[u_{i,j} \, \lvert \, i,j \in \Nat^+]].
\end{displaymath}

We say that \(R\) is 
\emph{\(\omega\)-Koszul} \footnote{
Recall that a commutative \(\setN\)-graded \(\setC\)-algebra \(R\) is 
    \emph{Koszul} if \(\Tor_R^i(\setC,\setC)_j = 0\) if \(i \neq j\).
    If \(R\) is Koszul, then \(P_R(t,u) = R(-tu)^{-1}\). An important
    result by Fröberg\cite{Froeberg:Poincare} is that if
    \(R=\frac{\setC[x_1,\dots,x_n]}{(m_1,\dots,m_r)}\), with \(m_i\)
      monomials (not of degree 1) 
      then \(R\) is Koszul if and only if all \(m_i\) are quadratic.}
if for every term \(t^i \vektor{u}^\vektor{\alpha}\) occurring in
\(P_R^M(t,\vektor{u})\), \(\tdeg{\vektor{u}}=i\).

\begin{lemma}\label{lemma:pbapprox}
  Let \(W \in \oksets\), and let
  \begin{displaymath}
    V_1 \subseteq V_2 \subseteq V_3 \subseteq \cdots \subseteq W, \quad
    \cup_{n=1}^\infty V_n = W
  \end{displaymath}
  be an ascending family in \(\oksets^f\).
  Then
  \begin{align}
    \label{eq:pblimit1}
      \lim_{n \to \infty}   \UNI_{V_n}(\vektor{u}) & =
      \UNI_W^f(\vektor{u})\\ 
    \label{eq:pblimit2}
      \lim_{n \to \infty}
      P_{\setC[Y(V_n)]}^{\UNI_{V_n}}(t,\vektor{u}) & = 
      P_{\UNI_{\setC[Y(W)]}}^{\UNI_W^f}(t,\vektor{u})\\
    \label{eq:pblimit3} 
      \lim_{n \to \infty}  P_{\UNI_V}(t,\vektor{u}) & =
      P_{\UNI_{W}^f}(t,\vektor{u})
  \end{align}
  where the topology on  \(\Z[[u_{i,j} \, \lvert \, i,j \in \Nat^+]]\) and
  on
  \(\Z[[t]][[u_{i,j} \, \lvert \, i,j \in \Nat^+]]\) is that of
  point-wise convergence.
\end{lemma}

\begin{definition}
    Let \(\vektor{x}=(x_1,x_2,x_3,\dots)\).
  For \(\sigma \subseteq \setN^+\) we put
  \begin{equation}
    \label{eq:sigm1}
    \vektor{x}_\sigma = \prod_{i \in \sigma} x_i,
  \end{equation}
  and if \(j \in  \sigma\), we 
  put 
  \begin{equation}
    \label{eq:sigm2}
  \vektor{x}_{\sigma,j} = x_j \vektor{x}_\sigma.  
  \end{equation}
\end{definition}

\begin{lemma}\label{lemma:nkoz}
  Let \(R=\setC[x_1,\dots,x_n]\), and let \(\mathbf{m}=(x_1,\dots,x_n)\)
  denote the graded maximal ideal. Put \(M=R/\mathbf{m}^2\). 
  Then 
  \begin{equation}
    \label{eq:sillyHilb}
    M(u_1,\dots,u_n) = 1+ u_1 + \dots + u_n
  \end{equation}
  \(M\) is Koszul, so
  \begin{equation}
    \label{eq:absolute}
    P_M(t,u_1,\dots,u_n)  = 
    \frac{1}{M(-tu_1,\dots,-tu_n)} 
      = \frac{1}{1-tu_1 - \cdots - tu_n}
  \end{equation}
  Furthermore,
  \begin{equation}
    \label{eq:relative}
      P_R^M(t,u_1,\dots,u_n)  = 1 + 
      \sum_{\substack{\sigma \subseteq  \set{1,\dots,n}\\
          1 \leq \tdeg{\sigma} \leq n\\
          j \in \sigma}} 
      t^{\tdeg{\sigma}} \vektor{u}_{\sigma,j} 
      + 
      \sum_{\substack{\sigma \subseteq  \set{1,\dots,n}\\
          2 \leq \tdeg{\sigma} \leq n}}
      (\tdeg{\sigma} - 1) t^{\tdeg{\sigma}-1} \vektor{u}_\sigma
  \end{equation}
\end{lemma}
\begin{proof}
  As a \(\C\)-vector space, \(M\) is spanned by \(1,x_1,\dots,x_n\),
  so \eqref{eq:sillyHilb} follows. Since the defining ideal of \(M\)
  is a quadratic monomial ideal, a theorem by Fröberg
  \cite{Froeberg:Poincare}  shows that \(M\) is Koszul. Hence
  \eqref{eq:absolute} follows.

  For the result on the relative Poincaré-Betti series, we'll use a
  result by Bayer \cite{MID}, which yields that the coefficient of the
  \(t^i\vektor{u}\)-term is equal to the rank of the torsion-free part
  of the \((i-2)\)'th reduced homology of the complex
  \begin{equation}
    \label{eq:complexBayer}
      \Delta_{\vektor{u}} = \setsuchas{F \subset
        [n]}{\vektor{x}^{\vektor{u} - F} \in \mathbf{m}}
  \end{equation}
  Here, \(F\) is identified with its characteristic vector, which is the
  square-free vector in \(\Nat^n\) with a \(1\) in position \(i\) if
  \(i \in F\), and zero otherwise.

  If \(\vektor{u}\) is square-free, \(\tdeg{\vektor{u}} = k\), then 
  \begin{displaymath}
    F \in \Delta_{\vektor{u}} \iff \# \left(\vektor{u} - F \right)
    \ge 2 \iff \#F \le \tdeg{\vektor{u}} - 2 = k-2,
  \end{displaymath}
  so \(\Delta_{\vektor{u}}\) is the \((k-3)\)-skeleton of a
  \((k-1)\)-simplex, hence the reduced homology is concentrated in
  degree \(k-3\), where it is \(\Z^{k-1}\).

  If \(\vektor{x}^{\vektor{u}} = x_a^2 \vektor{x}^{\vektor{v}}\),
  where \(\vektor{v}\) is square-free, and does not contain \(a\), and
  if \(\tdeg{\vektor{u}} = k\), (so \(\tdeg{v}=k-2\)),  then 
  \begin{displaymath}
    F \in \Delta_{\vektor{u}} \iff F \subseteq \vektor{v} \quad \text{
      or } 
    \quad 
    a \in F \text{ and } \# F \le k-3,
  \end{displaymath}
  thus \(\Delta_{\vektor{u}}\) is the disjoint union of a
  \((k-3)\)-simplex, which has no reduced homology, and a
  \((k-3)\)-sphere, which has homology concentrated in degree
  \((k-3)\), where it is \(\Z\).

  If \(\divides{x_a^2}{\vektor{x}^{\vektor{u}}}\) for two different
  \(a\)'s, so that e.g. \(\vektor{x}^{\vektor{u}} = x_1^2 x_2^2
  \vektor{x}^{\vektor{v}}\), then 
  \begin{displaymath}
    \forall F \subset [n]: \quad \vektor{x}^{\vektor{u} - F} \in
    \mathbf{m}, 
  \end{displaymath}
  thus
  \begin{displaymath}
    \Delta_{\vektor{u}} = 2^{[n]}
  \end{displaymath}
  is acyclic. 

  Summarizing we have that 
  \begin{equation}
    \label{eq:TorM}
    \dim_\C \Tor_{i,\vektor{u}}^R(M,\C) = 
    \begin{cases}
      i & \vektor{u} \text{ square-free of total degree } i+1 \\
      1 & u_a = 2, \, u_j \le 1 \text{ for } j \neq a \text{ and }
      \tdeg{\vektor{u}} = i+1 \\
      0 & \text{ otherwise}
    \end{cases}
  \end{equation}
  From this, \eqref{eq:relative} follows.
\end{proof}

\begin{theorem}\label{thm:UNIpoincare}
  \(\UNI^f\) is \(\omega\)-Koszul, and
  \begin{equation}
    \label{eq:kosz}
      P_{\UNI^f} = \frac{1}{\UNI^f(-t\vektor{u})}
      = \prod_{i=1}^\infty \frac{1}{1-t \sum_{j=1}^\infty u_{i,j}}
  \end{equation}

  Define
  \begin{equation}
    \label{eq:AB}
    A(t,\vektor{x}) = 
1 + 
      \sum_{\substack{\sigma \subseteq  \setN^+\\
          1 \leq \tdeg{\sigma} \\
          j \in \sigma}} 
      t^{\tdeg{\sigma}} \vektor{x}_{\sigma,j} 
      + 
      \sum_{\substack{\sigma \subseteq  \setN^+\\
          2 \leq \tdeg{\sigma} }}
      (\tdeg{\sigma} - 1) t^{\tdeg{\sigma}-1} \vektor{x}_\sigma
  \end{equation}
  For each positive integer \(i\), put 
  \[\vektor{u}_i =  (u_{i,1},u_{i,2}, u_{i,3},\dots).\] 
  Put 
  \begin{displaymath}
      \vektor{u}  =(\vektor{u}_1,\vektor{u}_2,\vektor{u}_3,\dots) 
      =(u_{1,1},u_{1,2},\dots; u_{2,1}, u_{2,2},
  \dots;u_{3,1}, u_{3,2},\dots; \dots)
  \end{displaymath}

  Then
  \begin{equation}
    \label{eq:poincareUNI}
    P_{\setC[Y^*]}^{\UNI^f}(t,\vektor{u})  = 
    \prod_{i=1}^\infty A(t,\vektor{u}_i)
  \end{equation}

  Similarly: \(\UNI_{\sqfree}^f\) is \(\omega\)-Koszul, and 
  \begin{equation}
    \label{eq:pbsqfree}
    P_{\UNI_{\sqfree}^f} = \frac{1}{\UNI_{\sqfree}^f(-t\vektor{u})}
    = \prod_{i=1}^\infty \frac{1}{1-t u_{i,1}} 
  \end{equation}
  \begin{equation}
    \label{eq:obsqfree2}
    P_{\setC[y_{1,1}, y_{2,1}, y_{3,1}, \dots]}^{\UNI_{\sqfree}^f}(t,\vektor{u})  = 
    A(t,\vektor{u}_1)
  \end{equation}
\end{theorem}
\begin{proof}
  We apply Theorem~\ref{lemma:pbapprox}, putting \(W=\setN^+\), \(V_n=\) the
  sub-monoid of \((\setN^+,\oplus)\) generated by \(\setsuchas{p_i^j}{1
    \leq i,j \leq n}\). Then 
  \begin{equation}\label{eq:square}
    \begin{split}
    \UNI_{V_n}^f &= \frac{\setC[\setsuchas{y_{i,j}}{1 \leq i,j \leq
        n}]}{(\setsuchas{y_{i,j} y_{i,j'}}{1 \leq i,j,j' \leq n})} \\
    &=       \frac{\setC[\setsuchas{y_{1,j}}{1 \leq j \leq
        n}]}{(\setsuchas{y_{1,j} y_{1,j'}}{1 \leq j,j' \leq n})}
    \otimes_\setC \cdots   \otimes_\setC 
     \frac{\setC[\setsuchas{y_{n,j}}{1 \leq j \leq
        n}]}{(\setsuchas{y_{n,j} y_{n,j'}}{1 \leq j,j' \leq n})} \\
    &= \left( \frac{\setC[x_1,\dots,x_n]}{\mathbf{m}^2} \right)^{\otimes n}
    \end{split}
  \end{equation}
  By Lemma~\ref{lemma:nkoz},
  \(\frac{\setC[x_1,\dots,x_n]}{\mathbf{m}^2}\) is Koszul and has
  absolute Poincaré-Betti series
  \begin{displaymath}
    \frac{1}{1-tu_1 - \cdots - tu_n},
  \end{displaymath}
  hence 
  \(\UNI_{V_n}^f\) is Koszul and 
  \begin{displaymath}
    P_{\UNI_{V_n}^f}(t,u_{1,1},\dots,u_{n,n}) = 
    \prod_{i=1}^n\frac{1}{1-tu_{i,1} - \cdots - tu_{i,n}}.
  \end{displaymath}
  It follows that
  \begin{displaymath}
    \begin{split}
    P_{\UNI^f}(t,\vektor{u}) & = 
      \lim_{n \to \infty} P_{\UNI_{V_n}^f}(t,u_{1,1},\dots,u_{n,n}) \\
      &= \lim_{n \to \infty} \prod_{i=1}^n\frac{1}{1-tu_{i,1} - \cdots
        - tu_{i,n}} \\
      &= \prod_{i=1}^\infty\frac{1}{1-t\sum_{j=1}^\infty u_{i,j}}.
    \end{split}
  \end{displaymath}
  The assertion about the relative Poincaré-Betti series is proved in
  a similar way. The square-free case is analogous.
\end{proof}

\begin{example}
  We continue our study of Example~\ref{example:pairwise}. Let \(V_n\) be
  as above, so that 
  \[\setC[Y(V_n)] = \setC[x_1,\dots,x_n]^{\otimes n},\]
  and \(\UNI_{V_n}^f\) is as in \eqref{eq:square}. The induced
  bi-grading on the \(i\)'th copy of \(C[x_1,\dots,x_n]\) gives
  \(x_j\) bi-grade \((1,p_i^j)\), hence 
  \begin{displaymath}
    \label{eq:Hilb}
    \UNI_{V_n}^f(u_1,u_2) = \prod_{i=1}^n \left(1 +
      \sum_{j=1}^n u_1u_2^{p_i^j} \right).
  \end{displaymath}

    We have that 
\begin{displaymath}
  \begin{split}
    \UNI_{V_1}^f(u_1,u_2) &= 1 + u_2^2 u_1\\
    \UNI_{V_2}^f(u_1,u_2) &=  1 + u_1u_2^2 
+ u_1u_2^3 
+ u_1u_2^4
+u_1^2u_2^5
+u_1^2u_2^7
+ u_1u_2^9
+ u_1^2u_2^{11}
+u_1^2u_2^{13}
 \\
    \UNI_{V_3}^f(u_1,u_2) &=  1 + u_1u_2^2 + u_1u_2^3 + u_1u_2^4 +
    (u_1+u_1^2)u_2^5 + \dots 
  \end{split}
\end{displaymath}
and \(\UNI_{V_n}^f(u_1,u_2)\) converges to \eqref{eq:gold} as \(n \to
\infty\). 

Now collapse this grading to an \(\Nat\)-grading.
Since \(\UNI_{V_n}\) is Koszul, we have that 
\begin{displaymath}
  P_{\UNI_{V_n}}(t,u) = \frac{1}{\UNI_{V_n}(-tu)} 
\end{displaymath}
which converges to 
\begin{displaymath}
  \begin{split}
  \frac{1}{\UNI^f(-tu)} &=  \frac{1}{1 -t^2u^2 - t^3u^3 -t^4u^4
  -2t^5u^5 -3t^7u^7 - \dots}\\ & = 1+u^2t^2+u^3t^3+ 2u^4t^4+ \dots
  \end{split}
\end{displaymath}
We have that
\begin{displaymath}
  \begin{split}
    P_{\UNI_{V_1}}(t,u) &=  \frac{1}{1- t^2u^2} = 1 + t^2u^2 + t^4u^4 +
    t^6u^6 + t^8u^8 + \dots \\
    P_{\UNI_{V_2}}(u) &=  \frac{1}{1- t^2u^2 - t^3u^3  - t^4u^4 - t^5u^5 -
      t^7u^7 - t^{9}u^{9} -  t^{11}u^{11}} t^{13}u^{13}  \\
    &= 1 + u^2t^2 + u^3t^3 + 2u^4t^4 + 3u^5t^5 + 4u^6t^6 + 8u^7t^7 +
    10u^8t^8+ \dots
  \end{split}
\end{displaymath}
Finally, 
\begin{displaymath}
  \begin{split}
  P_{\UNI^f}(t,u) &= \prod_{i=1}^\infty A(t, u^{p_i}) \\
  &= \prod_{i=1}^\infty 
\left( 1 + 
      \sum_{\substack{\sigma \subseteq  \setN^+\\
          1 \leq \tdeg{\sigma} \\
          j \in \sigma}} 
      t^{\tdeg{\sigma}} u^{p_i^j + \sum_{v \in \sigma} p_i^v}
      + 
      \sum_{\substack{\sigma \subseteq  \setN^+\\
          2 \leq \tdeg{\sigma} }}
      (\tdeg{\sigma} - 1) t^{\tdeg{\sigma}-1} u^{\sum_{v \in \sigma}
        p_i^v} 
    \right)
  \end{split}
\end{displaymath}
The reader is cordially invited to simplify this expression, and to
compute its first few terms.
\end{example}

\end{section}

\begin{section}{The case of a finite $V$}
  \begin{subsection}{Simple properties}
    We now assume that \(V \in \oksets^f\). To avoid trivial special
    cases, we assume that \(V\) contains at least two elements.
    Then \(\UNI_V = \UNI_V^f\)
    is a local Artin ring, where the maximal ideal is spanned (as a
    vector space) by \(\setsuchas{e_j}{j \in V \setminus
      \set{1}}\). It is easy to see that the maximal ideal is
    nilpotent, thus elements of \(\UNI_V\) are either units or
    nilpotent. 
\end{subsection}

  \begin{subsection}{Presentations}
  \begin{theorem}\label{thm:presentationsV}
    Let \(V \in \oksets^f\), \(\tdeg{V}=r\).
    \(\UNI_V\) is a cyclic
    module over \(\setC[Y(V)]\).
    Furthermore, if we define the following ideals of \(\setC[Y(V)]\), 
    \begin{equation}
      \label{eq:ids}
      \begin{split}
        A_V &= \setC[Y(V)] \setsuchas{y_{i,j}^2}{y_{i,j} \in Y(V)}\\
        B_V &= \setC[Y(V)] \setsuchas{y_{i,j} y_{i,j'}}{y_{i,j},
          y_{i,j'}\in Y(V), \, j < j'} \\
        C_V &= \setC[Y(V)]       \setsuchas{\prod_{\ell=1}^r y_{i_\ell,j_\ell}}
      {\prod_{\ell=1}^r p_{i_\ell}^{j_\ell} \not \in V, \,
        \forall 1 \leq v \leq r: \prod_{\substack{1 \leq \ell \leq r\\ \ell
          \neq j}} p_{i_\ell}^{j_\ell} \in V}
      \end{split}
    \end{equation}
    then the indicated generating sets are in fact the unique
    multi-graded minimal generating sets, and
    \begin{equation}
      \label{eq:pres}
      \UNI_V = \frac{C[Y(V)]}{A_V + B_V + C_V}
    \end{equation}
    If we let \(W\) denote the monoid ideal that \(V\) generates, then
    \begin{equation}
      \label{eq:Wpres}
    \UNI_W = \frac{C[Y(V)]}{A_V + B_V },
    \end{equation}
    and \(\UNI_V\) is a cyclic \(\UNI_W \)-module, since
    \[\UNI_V  = \frac{\UNI_W }{C_V}.\]
    We put 
    \begin{displaymath}
      \setC[\overline{Y(V)}] = \frac{\setC[Y(V)]}{A_V} \simeq 
    \left( \frac{\setC[x]}{(x^2)}\right)^{\otimes r},
    \end{displaymath}
    an \(r\)-fold tensor power of the ring of dual numbers.
Then
    \(\UNI_V\) is also a cyclic \(\setC[\overline{Y(V)}]\)-module, 
    since
    \(\UNI_V  = \frac{\setC[\overline{Y(V)}]}{B_V +C_V}\).
  \end{theorem}
  \begin{proof}
    We have that \(W \in \oksets^f \cap \upsets^f\), and
    \(Y(V)=Y(W)\), so
    by Theorem~\ref{thm:pres}, \(\UNI_V\) and \(\UNI_W\) are cyclic
    \(\C[Y(V)]\)-modules. Clearly, \(\UNI_W\) has a \(\C\)-basis
    consisting of 
    \begin{displaymath}
      \setsuchas{e_{k}}{k \in W} = \setsuchas{e_{i_1} \oplus e_{i_2}
        \oplus \cdots \oplus e_{i_r}}{i_1,\dots,i_r \in V \cap \PP}
    \end{displaymath}
    which in \(C[Y(W)]\) corresponds to separated monomials in the
    variables in \(Y(W)\). From this, \eqref{eq:Wpres} follows, since
    \(A(V)+B(V)\) is precisely what wee need to divide out with in
    order to have only separated monomials. Similarly, 
\(\UNI_V\) has a \(\C\)-basis
    consisting of 
    \begin{displaymath}
      \setsuchas{e_{k}}{k \in V} = \setsuchas{e_{i_1} \oplus e_{i_2}
        \oplus \cdots \oplus e_{i_r}}{i_1,\dots,i_r \in V \cap \PP, \,
        i_1 \oplus \cdots \oplus i_r \in V}
    \end{displaymath}
    so the defining ideal of \(\UNI_V\) in \(\C[Y(V)]\) contains, in
    addition, those separated monomials (in \(y_{i,j}\)-variables) which
    do not correspond to an element in  \(V\); from this observation,
    and Theorem~\ref{thm:mingens},
    \eqref{eq:pres} follows.

    The remaining assertions follow trivially.
  \end{proof}

    \begin{example}
Let 
\begin{displaymath}
V=[10]=\set{1,2,3,4,5,6,7,8,9,10}.
\end{displaymath}
Then
\begin{displaymath}
  Y([10]) =
\set{y_{1,1}, \, y_{1,2}, \, y_{1,3}, \, y_{2,1}, \, y_{2,2}, \,
  y_{3,1}, \, y_{4,1}},
\end{displaymath}
  and
  \begin{displaymath}
A_{[10]} = (
y_{1,1} y_{1,1}, \,
y_{1,2} y_{1,2}, \,
y_{1,3} y_{1,3}, \,
y_{2,1} y_{2,1}, \,
y_{2,2} y_{2,2}, \,
y_{3,1} y_{3,1}, \,
y_{4,1} y_{4,1}),
  \end{displaymath}
  \begin{displaymath}
B_{[10]}  =(
{y_{1,1}}^2, \,
{y_{1,1}}^2, \,
{y_{1,2}}^2,\,
{y_{2,1}}^2),
  \end{displaymath}
and finally 
\begin{multline*}
C_{[10]}  = (
y_{1,1} y_{2,2}, \,
y_{1,1} y_{4,1}, \,
y_{2,2} y_{3,1}, \,
y_{1,2} y_{2,1}, \,
y_{1,3} y_{2,1}, \,
y_{2,1} y_{3,1}, \,
y_{2,1} y_{4,1}, \\ 
y_{1,2} y_{2,2}, \,
y_{1,2} y_{3,1}, \,
y_{1,2} y_{4,1}, \,
y_{1,3} y_{3,1}, \,
y_{1,3} y_{4,1}, \,
y_{2,2} y_{4,1}, \,
y_{1,3} y_{2,2}, \,
y_{3,1} y_{4,1})
\end{multline*}
We have that 
\begin{displaymath}
  \UNI_{[10]} \simeq \frac{C[y_{1,1}, \, y_{1,2}, \, y_{1,3}, \, y_{2,1}, \, y_{2,2}, \,
  y_{3,1}, \, y_{4,1}]}{A_{[10]} + B_{[10]} + C_{[10]}}
\end{displaymath}
has a \(\C\)-basis consisting of (the images of) the following
separated monomials: 
\begin{displaymath}
1,\,  y_{1,1}, \, y_{1,2}, \, y_{1,3}, \, y_{2,1}, \, 
y_{1,1}y_{2,1}\,
y_{2,2}, \,
  y_{3,1}, \, y_{4,1}, \,
y_{1,1}y_{3,1}.
\end{displaymath}
\end{example}
    
  \end{subsection}

  \begin{subsection}{The associated simplicial  complex}
    For terminology and general results regarding simplicial complexes
    and Stanley-Reisner rings, see \cite{Stanley:CombCom}.

    \begin{definition}
    Let \(S\) be a finite set. A simplicial complex on \(S\) is a
    subset \(\Delta \subseteq 2^S\) of the power-set of \(S\), such that 
    \begin{enumerate}[(i)]
    \item \(\set{s} \in \Delta\) for all \(s \in S\),
    \item \(\tau \subseteq \sigma \in \Delta \implies \tau \in \Delta\).
    \end{enumerate}
    If \(U \subseteq S\), we denote by \(\Delta_U\) the simplicial
    complex on \(U\) given by \(\Delta_U = \Delta \cap 2^U\).

    If the reduced homology \(\widetilde{H}^i(\Delta,\setC)\) is non-zero 
    for \(i=j\), but vanishes for \(i>j\) we say that \(\Delta\) has
    homological degree \(i\).
  \end{definition}

  \begin{definition}
    Let \(S=\set{s_1,\dots,s_n}\) be a finite set, and let
    \(S_1,\dots,S_r\) be a 
    partition of \(S\), i.e. \(S=\sqcup_{i=1}^r S_i\), \(S_i \cap S_j =
    \emptyset\) if \(i \neq j\). A partitioned simplicial complex
    \(\Delta\) 
    (corresponding to the partition \(S=\sqcup_{i=1}^r S_i\)) on
    \(S\) is
    an order ideal in the sub-poset
    \begin{equation}
      \label{eq:partsimp}
      \overline{2^S} = \setsuchas{\sigma \subseteq S}{\forall 1 \leq i
        \leq r: \# (\sigma 
        \cap S_i) \leq 1} \subseteq 2^S
    \end{equation}
    such that \(\set{s} \in \Delta\) for every \(s
    \in S\).

    We define the following ideals in \(\setC[\vektor{x}] =
    \setC[x_1,\dots,x_n]\): 
    \begin{equation}
      \label{eq:idsX}
      \begin{split}
      A_S &= \setC[\vektor{x}] \setsuchas{x_a^2}{1 \leq a \leq n } \\
      B_S &= \setC[\vektor{x}] \setsuchas{x_ax_b}{\exists i: s_a,s_b \in
        X_i} \\ 
      C_\Delta &= \setC[\vektor{x}]  \setsuchas{x_{a_1} \cdots x_{a_v}}{1 \leq
        a_1 < a_2 < \cdots < a_v \leq n, \set{a_1,\dots,a_v} \in
        \overline{2^X} \setminus \Delta}
      \end{split}
    \end{equation}
    We define the Stanley-Reisner algebra of \(\Delta\) as
    \begin{equation}
      \label{eq:SR}
          \setC[\Delta]=\frac{\setC[x_1,\dots,x_n]}{B_S + C_\Delta} 
    \end{equation}
    and the   Artinified Stanley-Reisner ring as 
\begin{equation}
      \label{eq:ASR}
          \setC[\overline{\Delta}]=\frac{\setC[x_1,\dots,x_n]}{A_S +B_S +
            C_\Delta}  
    \end{equation}
  \end{definition}

Note that a partitioned simplicial complex is also a simplicial
complex in the ordinary sense, that is, an order ideal in \(2^X\), and
that 
\(\setC[\Delta]\) coincides with the usual definition; the only
difference is that we have partitioned the Stanley-Reisner ideal
\(I_\Delta = B_X + C_\Delta\). The Artinified Stanley-Reisner ring
\[\setC[\overline{\Delta}] \simeq \frac{\setC[\vektor{x}]}{A_X + I_\Delta}\]
is identical to the one introduced by Sköldberg
\cite{Skold:Golod}. The only difference is that we may choose to
regard it as a cyclic \(\frac{\setC[\vektor{x}]}{A_X + B_X}\)-module as
well as 
a cyclic \(\setC[\vektor{x}]\)-module.

\begin{example}
  \begin{displaymath}
    \Delta([10]) = \set{\emptyset, \set{2}, \set{3}, \set{4}, \set{5},
      \set{2,3}, \set{7}, \set{8}, \set{9}, \set{2,5}}
  \end{displaymath}
looks like Figure~\ref{fig:Delta10}. 
\begin{figure}[htbp]
    \begin{center}
    \setlength{\unitlength}{1truecm}
      \begin{picture}(4,4)
        \put(0,1){\circle*{0.1}}
        \put(1,1){\circle*{0.1}}
        \put(3,1){\circle*{0.1}}
        \put(4,1){\circle*{0.1}}
        \put(0,0.5){4}
        \put(1,0.5){7}
        \put(3,0.5){8}
        \put(4,0.5){9}

        \put(2,1){\circle*{0.1}}
        \put(0,3){\circle*{0.1}}
        \put(4,3){\circle*{0.1}}
        \put(2,0.5){2}
        \put(0.1,3.1){3}
        \put(3.6,3.1){5}

        \put(2,1){\line(-1,1){2}}
        \put(2,1){\line(1,1){2}}
        \put(1,2.2){6}
        \put(2.4,2.2){10}
      \end{picture}
      \caption{\(\Delta([10])\)}
      \label{fig:Delta10}
    \end{center}
  \end{figure}

\end{example}

For subsequent use, we state the following theorem:

\begin{theorem}\label{thm:bettiASR} 
  Let \(\Delta\) be a simplicial complex, and 
  \(\setC[\overline{\Delta}]\) its Artinified Stanley-Reisner ring.
Let \(\C[\overline{\vektor{x}}]\) denote the (smallest) polynomial
ring of which 
\(\setC[\overline{\Delta}]\)  is an epimorphic image.
  Then the multi-graded Betti numbers 
  \begin{displaymath}
  \beta_{i, \vektor{\alpha}} = \dim_\C
  \Tor_{\vektor{\alpha}}^{\C[\overline{\vektor{x}}]}
  (\setC[\overline{\Delta}], \C)
  \end{displaymath}
  are given by
  \begin{equation}
    \label{eq:1}
      \beta_{i,\vektor{\alpha}} =  
\dim_\setC     \widetilde{H}^{\tdeg{\vektor{\alpha}}-i-1}(\Delta_{U};
\, \setC), \quad 
\text{ where }  U=\supp(\vektor{\alpha})
  \end{equation}
\end{theorem}
  \begin{proof}
    This is the formula for the Betti numbers of the corresponding
    indicator algebra \(\C\{\Delta\}\) \cite{Aramova:Cohomology,
      Aramova:Gotzman}. This skew-commutative algebra has the same
    Betti numbers as \(\setC[\overline{\Delta}]\) \cite{Skold:Golod}.
  \end{proof}

\begin{theorem}
  For \(V \in \oksets^f\), let \(P=V \cap \PP\), 
\(Q = V \cap \Primes\).
Number the elements in \(Q\) as \(q_1,\dots,q_r\) such
  that \(q_1 < q_2 < \cdots < q_r\).
  Partition \(P=P_1 \cup P_2 \cdots \cup P_r\), \(P_i =
  \setsuchas{q_i^j}{q_i^j \in P}\). Define a partitioned simplicial
    complex \(\Delta(V)\) on \(P\) by 
    \begin{equation}
      \label{eq:simpcomp}
      \set{q_{i_1}^{a_1},\dots,q_{i_s}^{a_s}} \in \Delta(V) \quad \iff
      \quad q_{i_1}^{a_1} \cdots q_{i_s}^{a_s} \in V
    \end{equation}
    Then 
    \begin{equation}
      \label{eq:simisaro}
      \setC[\overline{\Delta(V)}] \simeq \UNI_V 
    \end{equation}
    as graded \(\setC\)-algebras,
    and
    \begin{equation}
      \label{eq:posi}
      (\Delta(V), \subseteq)  \simeq (V, \leq_{\oplus_V})      
    \end{equation}
    as posets.
\end{theorem}
\begin{proof}
  Since \eqref{eq:simpcomp} gives a bijection between the basis vectors
  of \(\setC[\overline{\Delta(V)}]\) and those of \( \UNI_V\), we need
  only check that the multiplication is the same. 
  Let \(\sigma, \tau \in \Delta(V)\),
and define
  \begin{displaymath}
    \begin{split}
    \sigma(i) &= 
    \begin{cases}
      1 & q_i \in \sigma \\
      0 & q_i \not \in \sigma
    \end{cases}
\\
      x_\sigma & = \prod_{i=i}^r x_i^{\sigma(i)}\\
      e_\sigma &=  e_m \quad \text{ where } m= {\prod_{j \in \sigma} j}
    \end{split}
  \end{displaymath}
  and similarly for \(\tau\).
  Then 
  \begin{displaymath}
    e_\sigma e_\tau = 
    \begin{cases}
      e_{\sigma \cup \tau} & \text{ if } \sigma \cap \tau = \emptyset,
      \, \sigma \cup \tau \in \Delta(V) \\
      0 & \text{ otherwise}
    \end{cases}
  \end{displaymath}
  Furthermore \(x_\sigma, x_\tau \in \setC[\overline{\Delta(V)}]\), 
and
  \begin{displaymath}
    x_\sigma x_\tau = 
    \begin{cases}
      x_{\sigma \cup \tau} & \text{ if } \sigma \cap \tau = \emptyset,
      \, \sigma \cup \tau \in \Delta(V) \\
      0 & \text{ otherwise}
    \end{cases}
  \end{displaymath}
  This shows that the  multiplications are the same.

  The bijection \eqref{eq:simpcomp} is easily seen to be the desired
  poset isomorphism \eqref{eq:posi}.
\end{proof}

\end{subsection}

\begin{subsection}{The socle}
There is a short exact sequence of complex vector spaces and linear maps
\begin{equation}
  \label{eq:monsyz}
  \begin{split}
  0 \to K_1(V) \to \UNI_V \otimes \UNI_V  & \to \UNI_V \to 0 \\
  f \otimes g & \mapsto f \oplus g
  \end{split}
\end{equation}
which restrict to 
\begin{equation}
  \label{eq:monsyz2}
  0 \to K_2(V) \to \UNI_V^+ \otimes \UNI_V^+   \to \UNI_V^+
\end{equation}

We denote by \(\pr\) the projection \(\UNI_V^+ \otimes \UNI_V^+   \to \UNI_V^+\)
to the first factor.

We call elements in \(K_1(V)\) of the form \(e_a \otimes e_b\), \(e_a
\oplus e_b = 0\),
\emph{monomial multiplicative syzygies} and those of the form 
\(e_a \otimes e_b - e_c \otimes e_d\) \emph{binomial multiplicative
  sysygies}. 

\begin{lemma}\label{lemma:monsyz}
  \(K_1(V)\) is spanned (as a \(\C\)-vector space) by monomial and
  binomial multiplicative syzygies; consequently, so is \(K_2(V)\).
  Let \(n=\tdeg{V} = \dim_\C \UNI_V\), then 
   \(\dim_\C K_1(V) = n^2 - n\), \(\dim_\C K_2(V) = (n-1)^2 - (n-1)\).
\end{lemma}
\begin{proof}
  Since the map \eqref{eq:monsyz} is multi-homogeneous, it is enough
  to study it in one multi-degree 
  \[\vektor{\alpha} = (\alpha_1,\dots,\alpha_r),\] where 
  \[ k = p_1^{\alpha_1} \cdots  p_r^{\alpha_r}  \in  V \oplus V.\]
  If 
  \begin{displaymath}
    f = \sum_{\substack{a,b \in V\\ab=k}} c_{a,b} e_a \otimes e_b 
  \end{displaymath}
  is an element of \(K_1(V)\) of multi-degree \(\vektor{\alpha}\),
  then we can write \(f=f_1 + f_2\) with 
  \begin{displaymath}
    f_1 = \sum_{\substack{a,b \in V\\ab=k\\ a \oplus_V b = 0}} c_{a,b} e_a
    \otimes e_b, \qquad
    f_2 = \sum_{\substack{a,b \in V\\ab=k\\ a \oplus_V b = k}} c_{a,b} e_a
    \otimes e_b.
  \end{displaymath}
  Then \(f_1\) is a linear combination of monomial syzygies, so it is
  enough to show that \(f_2\) is a linear combination of binomial
  syzygies. By construction, 
  \begin{displaymath}
    \sum_{\substack{a,b\\ab=k\\ a \oplus_V b = k}} c_{a,b} = 0.
  \end{displaymath}
  If we order the pairs \((a,b)\) with \(a,b \in V\), \(a \oplus_V b
  = k\) linearly, this becomes
  \begin{displaymath}
    \sum_{i=1}^L c_i = 0,
  \end{displaymath}
  or 
  \begin{displaymath}
    \left[
    \begin{matrix}
      1 & 1 & \cdots & 1
    \end{matrix}
    \right]
    \left[
      \begin{matrix}
        c_1 \\ c_2 \\ \vdots \\ c_L
      \end{matrix}
      \right]
      = 0
  \end{displaymath}
  It is an elementary fact that the solution set is spanned by 
  \begin{displaymath}
    \left\{
    \left[
      \begin{matrix}
        -1 \\ 1 \\ 0 \\ 0 \\ 0 \\ \vdots \\ 0
      \end{matrix}
      \right],
    \left[
      \begin{matrix}
        -1 \\ 0 \\ 1 \\ 0 \\ 0 \\ \vdots \\ 0
      \end{matrix}
      \right],
    \left[
      \begin{matrix}
        -1 \\ 0 \\ 0 \\ 1 \\ 0 \\ \vdots \\ 0
      \end{matrix}
      \right], \dots,
    \left[
      \begin{matrix}
        -1 \\ 0 \\ 0 \\0 \\  \vdots \\ 0 \\ 1
      \end{matrix}
      \right]
      \right\}
  \end{displaymath}
  Thus \(f_2\) can be written as a linear combination of binomial
  syzygies.

  Since dimension of vector spaces is additive function over short
  exact sequences, the assertions on the dimension of \(K_1\) and
  \(K_2\) follows from the fact that \(\dim_\C \UNI_V =n\), \(\dim_C
  \UNI_V \otimes \UNI_V = n^2\), \(\dim_\C \UNI_V^+ =n-1\), \(\dim_C
  \UNI_V^+ \otimes \UNI_V^+ = (n-1)^2\).
\end{proof}

\begin{definition}\label{def:socle}
The \emph{socle} of \(\UNI_V\) consists of all
elements in the maximal ideal that annihilates the maximal ideal. 
We denote it by \(\mathrm{Socle}(\UNI_V)\).
\end{definition}

\begin{lemma}\label{lemma:socle}
  The socle of \(\UNI_V\) is generated (as a \(\setC\)-vector space) by 
  \begin{multline}
    \label{eq:gensocle}
    \setsuchas{e_k, \, k \in V \setminus \set{1}}{km \not \in V \text{
        if } m \in V \setminus \set{1} \text{ and } \gcd(k,m) =1} = \\
    =
    \setsuchas{e_k}{e_k \otimes e_m \in K_2(V) \text{ for all } m \in \UNI_V^+}
  \end{multline}
  In fact, the above set is a basis.
\end{lemma}
\begin{proof}
  \(\UNI_V\) and its maximal ideal are multi-graded:
  \(\UNI_V\) has the \(\C\)-basis 
  \begin{displaymath}
    \setsuchas{e_k}{k \in V},
  \end{displaymath}
  and 
  \begin{displaymath}
    \setsuchas{e_k}{k \in V \setminus \set{1}}
  \end{displaymath}
  is a basis for the maximal ideal. Hence, the set 
  \begin{multline*}
    \setsuchas{e_k}{k \in V \setminus \set{1}, \, e_k \oplus_V e_m = 0
      \text{ for } m \in V \setminus \set{1}} = \\ =
    \setsuchas{e_k, \, k \in V \setminus \set{1}}{km \not \in V \text{
        if } m \in V \setminus \set{1} \text{ and } \gcd(k,m) =1} 
  \end{multline*}
  is a basis for the socle. This is precisely the set of \(e_k\), \(k
  \in V \setminus \set{1}\) such that \(e_k \otimes e_m\) is a
  multiplicative syzygy for all \(e_m \in \UNI_V^+\).
\end{proof}

  \begin{lemma}\label{lemma:soclefacet}
    Let \(k \in V\). The following are equivalent: 
    \begin{enumerate}[(i)]
    \item \(e_k \in \mathrm{Socle}(\UNI_V)\).
    \item \(k\) is a  maximal element in \((V, \leq_{\oplus_V})\).
    \item \(k\) corresponds to a facet in \(\Delta(V)\).
    \item \(\pr^{-1}(e_k) \subseteq K_2\).
    \end{enumerate}
  \end{lemma}
  \begin{proof}
    If \(k\) is maximal, then \(k \oplus_V m = 0\) for all \(m \in V
    \setminus \set{1}\), hence \(e_k  \oplus_V e_m = 0\) for all \(m \in V
    \setminus \set{1}\), hence \(e_k \in  \mathrm{Socle}(\UNI_v)\).
    If \(k\) is not maximal, so that \(k <_{\oplus_V} m\), then 
    \(m =   k \oplus_v c\) for some \(c \in V \setminus \set{1}\), 
    hence \(e_m =   e_k \oplus_v e_c\), hence \(e_k\) does not
    annihilate all elements in the maximal ideal, hence \(e_k \not \in
    \mathrm{Socle}(\UNI_v)\). 
    
    By \eqref{eq:posi},  maximal elements in \((V, \leq_{\oplus_V})\)
    correspond to facets in \(\Delta(V)\).
  \end{proof}

  Since \(\UNI_V\) is Artin, it is Cohen-Macaulay. It is Gorenstein
  if and only if the socle is one-dimensional \cite[Theorem
  12.4]{Stanley:CombCom}, hence 

  \begin{lemma}\label{lemma:whenGorenstein}
    If \(k=\max V\) then \(e_k \in \mathrm{Socle}(\UNI_V)\).
    The following are equivalent: 
    \begin{enumerate}[(i)]
    \item \(\UNI_V\) is Gorenstein,
    \item \(\mathrm{Socle}(\UNI_V) = \setC e_k\),
    \item For every \(j \in V \setminus \set{k,1}\) there exists at
      least one \(i=i(j)\in V \setminus \set{k,1}\) such that
      \(\gcd(i,j)=1\), \(ij \in V\).
    \item \(V\) is a principal order ideal.
    \item \(\Delta(V)\) is a simplex.
    \end{enumerate}
  \end{lemma}
  \begin{proof}
    If \(k=\max V\) then it is a maximal element in \((V,
    \geq_{\oplus_V})\), hence \(e_k \in \mathrm{Socle}(\UNI_V)\) by the
    previous Lemma. In fact, the only \(V \in \oksets^f\) which have only one
    maximal element are the principal order ideals.
  \end{proof}

  \end{subsection}

  \begin{subsection}{The Koszul property}
    \begin{theorem}
      Let \(V \in \oksets^f\) contain at least two elements,
      and let \(Q = V \cap \PP\). Let \(W\) denote the sub-monoid of
      \((\setN^+,\oplus)\) generated by \(Q\).
      The following are equivalent:
      \begin{enumerate}[(i)]
      \item \label{en:kos}
        \(\UNI_V\) is Koszul,
      \item \label{en:cV}
        \(C_V\) is 0 or generated in degree 2.
      \item \label{en:mg}
        If \(w \in W \setminus V\) then 
        \begin{displaymath}
 \exists q_1,q_2 \in  W: \quad
 \gcd(q_1,q_2)=1,  \quad
 \ddivides{q_1q_2}{w}, \quad
 q_1q_2 \in W \setminus V.          
        \end{displaymath}
      \item \label{en:betti1}
        For any \(U \subseteq Q\) with more than two elements,
        \[\tilde{H}^{\tdeg{U}-2}(\Delta(V)_U,\setC) = 0.\]
      \end{enumerate}
    \end{theorem}
    \begin{proof}
      By \eqref{eq:pres} and the fact that \(A_V\) and \(B_V\) are
      quadratic, Fröberg's result \cite{Froeberg:Poincare} that a
      monomial algebra is Koszul if and only if it is quadratic gives the
      equivalence between (\ref{en:kos}) and (\ref{en:cV}).
      (\ref{en:mg}) says that any  \(e_w\) in the defining ideal of
      \(\UNI_V\) as a cyclic \(\UNI_W\)-module is divisible by some
      \(e_{q_1q_2}\), i.e. by some quadratic element. This is
      equivalent to all minimal generators being quadratic.

      Finally, let \(\beta_{1,j}\) denote the number of minimal
      generators of \(A_V + B_V + C_V\) of degree \(j\). 
      Since \(A_v + B_V\) are generated in degree 2, it follows
      from Theorem~\ref{thm:bettiASR} that for \(j>2\),
      \begin{displaymath}
        \beta_{1,j} = \sum_{\substack{U \subseteq Q\\ \tdeg{U}=j}} \dim_{\setC} 
        \tilde{H}^{\tdeg{U}-2}(\Delta(V)_U,\setC).
      \end{displaymath}
      We want \(\beta_{1,j}=0\) for \(j>2\), so (\ref{en:betti1}) follows.
    \end{proof}

    \begin{corr}
      If  \(\UNI_V\) is Koszul  then  for \(d > 1\), \(\Delta(V)\) can
      not contain a 
        ``punctured \(d\)-simplex'' i.e. a subset of the form \(2^U
        \setminus U\) with \(\tdeg{U}=d+1\).
    \end{corr}
    \begin{proof}
      If it does, \(\Delta(V)_U \simeq \mathbb{S}^{d-1}\) so
      \(\tilde{H}^{d-1}(\Delta(V)_U,\setC) \neq 0\).
    \end{proof}

    \begin{corr}
      Let \(Q \subseteq \PP\) be finite and non-empty, and let \(V\) be
      the monoid ideal in \((\setN^+,\oplus)\) generated by \(Q\). Then
      \(V \in \oksets^f\) and \(\UNI_V\) is Koszul.
    \end{corr}
  \end{subsection}

  \begin{subsection}{Universality}
  \begin{theorem}\label{thm:everything}
  Let \(\Gamma\) be a finite simplicial complex. Then there exists
  infinitely many \(V \in \oksets\) such that \(\Delta(V) \simeq
  \Gamma\). 
\end{theorem}
\begin{proof}
  Without loss of generality we may assume that \(\Gamma\) is a
  simplicial complex on the set \(\set{1,\dots,n}\). Let
  \(q_1,\dots,q_n\) be any set of distinct primes, and define 
  \begin{displaymath}
    V = \setsuchas{\prod_{i \in \sigma} q_i}{\sigma \in \Gamma}.
  \end{displaymath}
  Then \(\Delta(V) \simeq \Gamma\).
\end{proof}

From this result we conclude that it isn't feasible to study the
(homological) properties of all truncations \(\UNI_V\), for \(V\)
finite. In a sequel to this article,  we'll restrict ourselves to the
special cases 
\(V=[n]\) and \(V=[n] \cap \sqfree\).
  \end{subsection}
\end{section}


\bibliographystyle{plain}
\bibliography{journals,articles,snellman}

\begin{thebibliography}{10}

\bibitem{Aramova:Cohomology}
Anetta Aramova, Luchezar~L. Avramov, and J{\"u}rgen Herzog.
\newblock Resolutions of monomial ideals and cohomology over exterior algebras.
\newblock {\em Transactions of the {A}merican {M}athematical {S}ociety},
  352(2):579--594, 1999.

\bibitem{Aramova:Gotzman}
Annetta Aramova, J{\"u}rgen Herzog, and Takayuki Hibi.
\newblock Gotzman theorems for exterior algebras and combinatorics.
\newblock {\em Journal of {A}lgebra}, 191:174--211, 1997.

\bibitem{MID}
Dave Bayer.
\newblock Monomial {I}deals and {D}uality.
\newblock Preprint, 1996.

\bibitem{Bourbaki:Commutative}
Nicolas Bourbaki.
\newblock {\em Commutative {A}lgebra}.
\newblock Hermann, 1972.

\bibitem{Cashwell:FPS}
E.~D. Cashwell and C.~J. Everett.
\newblock Formal power series.
\newblock {\em Pacific {J}ournal of {M}athematics}, 13:45--64, 1963.

\bibitem{Eliahou:MinRes}
S.~Eliahou and M.~Kervaire.
\newblock Minimal resolutions of some monomial ideals.
\newblock {\em Journal of {A}lgebra}, 129:1--25, 1990.

\bibitem{Froeberg:Poincare}
Ralf Fr{\"o}berg.
\newblock Determination of a class of {P}oincar\'e series.
\newblock {\em Math. Scand.}, 37(1):29--39, 1975.

\bibitem{Rib:GP}
Paulo Ribenboim.
\newblock Generalized power series rings.
\newblock In J.~Almeida, G.~Bordalo, and P.~Dwinger, editors, {\em Lattices,
  semigroups and universal algebra}, pages 15--33. Plenum, {N}ew {Y}ork, 1990.

\bibitem{Rib:UniZ}
Paulo Ribenboim.
\newblock Rings of generalized power series {II}: units and zero-divisors.
\newblock {\em Journal of {A}lgebra}, 168(1):71--89, Aug 1994.

\bibitem{dvar}
Emil~D. Schwab and Gheorghe Silberberg.
\newblock A note on some discrete valuation rings of arithmetical functions.
\newblock {\em Archivum Mathematicum (Brno)}, 36:103--109, 2000.

\bibitem{vring}
Emil~D. Schwab and Gheorghe Silberberg.
\newblock The valuated ring of the arithmetical functions as a power series
  ring.
\newblock {\em Archivum Mathematicum (Brno)}, 37(1):77--80, 2001.

\bibitem{Skold:Golod}
Emil Sk{\"o}ldberg.
\newblock Monomial {G}olod {Q}uotients of {E}xterior {A}lgebras.
\newblock {\em Journal of {A}lgebra}, 218:183--189, 1998.

\bibitem{Snellman:Dirichlet}
Jan Snellman.
\newblock Truncations of the ring of number-theorethic functions.
\newblock {\em Homology, {H}omotopy and {A}pplications}, 2:17--27, 2000.

\bibitem{Snellman:UniDivTop}
Jan Snellman.
\newblock The ring of arithmetical functions with unitary convolution:
  Divisorial and topological properties.
\newblock \emph{Research {R}eports in {M}athematics} 4/2002, {D}epartment of
  {M}athematics, {S}tockholm {U}niversity, apr 2002.

\bibitem{Stanley:CombCom}
Richard~P. Stanley.
\newblock {\em Combinatorics and {C}ommutative {A}lgebra}, volume~41 of {\em
  Progress in {M}athematics}.
\newblock Birkh{\"a}user, 2 edition, 1996.

\end{thebibliography}
\end{document}